\documentclass[11pt,a4paper]{article}
\usepackage[margin=25mm]{geometry}
\usepackage[T1]{fontenc}
\usepackage{lmodern,amsmath,amssymb,amsthm,mathtools,microtype}
\usepackage{xcolor,enumitem,fancyhdr}
\usepackage[nocompress]{cite}
\usepackage[colorlinks=true,linkcolor=blue!45!black,citecolor=blue!45!black,urlcolor=blue!45!black]{hyperref}
\hypersetup{pdftitle={A Liouville theorem for plurisubharmonic functions with periodic complex Monge-Ampere measure},pdfauthor={Jingchen Hu and Xushan Tu},pdfsubject={Periodic measures with locally Holder continuous Monge-Ampere subsolutions}}
\allowdisplaybreaks[2]
\numberwithin{equation}{section}
\newtheorem{theorem}{Theorem}[section]
\newtheorem{proposition}[theorem]{Proposition}
\newtheorem{lemma}[theorem]{Lemma}
\newtheorem{corollary}[theorem]{Corollary}
\theoremstyle{definition}
\newtheorem{definition}[theorem]{Definition}

\newcommand{\C}{\mathbb C}
\newcommand{\R}{\mathbb R}
\newcommand{\Hess}{D^2_{\C}}
\newcommand{\HR}{D^2_{\R}}
\newcommand{\tr}{\operatorname{tr}}
\newcommand{\osc}{\operatorname{osc}}
\newcommand{\cof}{\operatorname{cof}}
\newcommand{\Id}{I}
\newcommand{\G}{\mathcal G}
\newcommand{\B}{\mathcal B}
\newcommand{\Sph}{\mathcal S}
\newcommand{\dd}{i\partial\bar\partial}
\newcommand{\norm}[1]{\left\|#1\right\|}
\newcommand{\fint}{\mathop{\rlap{\,--}\!\int}\nolimits}
\newcommand{\PSH}{\operatorname{PSH}}
\newcommand{\MA}{\mathcal M}

\newcommand{\PolishOld}[1]{\begingroup\color{red}%
  \hypersetup{linkcolor=red,citecolor=red,urlcolor=red}%
  \if\relax\detokenize{#1}\relax\textit{[No corresponding text.]}\else#1\fi\endgroup}
\newcommand{\PolishNew}[1]{\begingroup\color{blue}%
  \hypersetup{linkcolor=blue,citecolor=blue,urlcolor=blue}%
  \if\relax\detokenize{#1}\relax\textit{[Deleted.]}\else#1\fi\endgroup}

\DeclareRobustCommand{\PolishMath}[3]{%
  \text{\textsuperscript{\scriptsize[#1]}}%
  {\color{red}\if\relax\detokenize{#2}\relax\text{\scriptsize[none]}\else#2\fi}%
  \,\to\,{\color{blue}#3}}

\title{\vspace{-1.2em}A Liouville theorem for plurisubharmonic functions\\
with periodic complex Monge--Amp\`ere measure}
\author{Jingchen Hu \qquad Xushan Tu}
\date{October 29, 2026}
\begin{document}
\maketitle
\vspace{-2em}
\begin{abstract}
\begingroup
We study entire continuous plurisubharmonic solutions of the complex Monge--Amp\`ere equation
\[
 (i\partial\bar\partial u)^n=2^n n!\,\mu\quad\text{in }\C^n,
\]
in the Bedford--Taylor sense, where $\mu\not\equiv0$ is a nonnegative locally finite Borel measure periodic with respect to a lattice of full real rank. We assume that $\mu$ is locally dominated by complex Monge--Amp\`ere measures of H\"older continuous plurisubharmonic functions, a condition under which the existence theorem of Demailly et al.~\cite{DDGHKZ} yields H\"older continuous periodic correctors. Under this assumption, we prove that every such solution with positive two-sided quadratic growth can be written uniquely, up to an additive constant, as the sum of a real quadratic polynomial and a H\"older continuous periodic function.
\par\endgroup
\end{abstract}

\begingroup\small
\par\endgroup

\section{Introduction and main results}

\subsection{The classification problem and related work}

The classical theorem of J\"orgens, Calabi, and Pogorelov~\cite{Jorgens,Calabi,Pogorelov} states that every entire classical convex solution of the real Monge--Amp\`ere equation with positive constant right-hand side is a quadratic polynomial. Related Liouville theorems for the real $k$-Hessian equations with positive constant right-hand side, of orders between two and one less than the dimension, were proved under positive quadratic lower growth, together with strict convexity or $(k+1)$-convexity~\cite{BCGJ2003,LRW2016}. More recently, Zhou and Zhu~\cite{ZZ2026} proved that every convex $C^2$ entire solution of these equations of order at least three is a quadratic polynomial, without any growth assumption.

For the complex Monge--Amp\`ere equation with positive constant right-hand side, the addition of an entire pluriharmonic function preserves both plurisubharmonicity and the equation. Further hypotheses are therefore needed to obtain quadraticity of the potential. Riebesehl and Schulz~\cite{RS1984} proved that the complex Hessian is constant when the associated K\"ahler metric is uniformly equivalent to the Euclidean metric. Wang~\cite{Wang2013} proved quadratic rigidity for entire plurisubharmonic solutions whose difference from a prescribed positive quadratic solution is of smaller than quadratic order at infinity. Li and Sheng~\cite{LiSheng2021} obtained a Liouville theorem under completeness of the associated K\"ahler metric and positive two-sided quadratic growth.

Recent results remove the completeness assumption under other natural hypotheses. Chen, Hu, and Sheng~\cite{CHS2026} proved that every smooth entire real convex solution is quadratic, without a growth assumption. Chu and Liu~\cite{ChuLiu2026} proved quadratic rigidity for smooth plurisubharmonic solutions under positive two-sided quadratic growth. Hu and Sheng~\cite{HuSheng2026} established this conclusion for locally bounded Bedford--Taylor solutions under the same growth condition and, more generally, under quadratic bounds with logarithmic factors. Their result gives the constant-density classification in the weak solution setting relevant to this paper.

When the coefficients or data of an equation are periodic, a natural question is how this periodicity is reflected in the structure of entire solutions satisfying suitable growth conditions.
Avellaneda and Lin~\cite{AL1989} and Li and Wang~\cite{LW2001} described solutions of polynomial growth to linear uniformly elliptic equations with periodic coefficients in terms of polynomials and periodic functions. For a class of periodic quasilinear elliptic equations, Moser and Struwe~\cite{MS1992} proved that every entire solution with asymptotically linear growth is an affine function plus a periodic function. Li and Liang~\cite{LiLiang2026} obtained a quadratic--periodic decomposition for a class of fully nonlinear elliptic equations.

For the real Monge--Amp\`ere equation, Y.~Y.~Li~\cite{Li1990} established existence and uniqueness results for periodic solutions. 
 Caffarelli and Li~\cite{CL2004} proved that every entire convex solution with a positive H\"older continuous periodic density is the sum of a quadratic polynomial and a periodic function. Li and Lu~\cite{LiLu2022} extended this classification to a class of nonsmooth positive periodic densities. Teixeira and Zhang~\cite{TZ2016} studied the asymptotic behavior of entire convex solutions with asymptotically periodic data. Zhang and Bao~\cite{ZB2018} obtained a parabolic analogue for ancient solutions with periodic spatial data. Wang~\cite{Wang2026} studied a coupled-drift Monge--Amp\`ere equation with periodic forcing.
Recently, Jin, Li, Tran, and Tu~\cite{JLT} established the classification for arbitrary nonzero nonnegative locally finite periodic Borel measures.

\subsection{Main results}

 We consider continuous entire plurisubharmonic (psh) solutions of
\begin{equation}\label{eq:normalization}
  \MA(u)=\mu\quad\text{in }\C^n, \qquad \MA(v):=\frac{(i\partial\bar\partial v)^n}{2^n n!},
\end{equation}
where $\MA(v)$ is defined for locally bounded psh functions $v$ in the Bedford--Taylor sense. For smooth $v$, we write $\Hess v=(v_{z_j\bar z_k})_{j,k=1}^n$, with $\partial_{z_j}=(\partial_{x_j}-i\partial_{y_j})/2$; thus $\MA(v)=\det(\Hess v)\,dV$.

Yau~\cite{Yau1978} established existence and uniqueness on compact K\"ahler manifolds for smooth positive data satisfying the mass compatibility condition. Ko{\l}odziej~\cite{Kol1998} proved existence of continuous potentials for nonnegative $L^p$ densities with $p>1$ satisfying the same compatibility condition, and later established their H\"older continuity~\cite{Kol2008}. Demailly, Dinew, Guedj, Hiep, Ko{\l}odziej, and Zeriahi~\cite[Theorem 4.3]{DDGHKZ} characterized measures admitting H\"older continuous potentials by local domination by Monge--Amp\`ere measures of H\"older continuous psh functions. In the present setting, this condition is
\begin{equation*}\tag{H}\label{eq:holder-measure}
 \begin{gathered}
 \text{For every }x\in\C^n,\text{ there exist a neighborhood }U\ni x\\
 \text{and a }v\in\PSH(U)\cap C^{0,\alpha}(U),
 \text{  }\alpha>0,\text{ such that }\mu|_U\le (i\partial\bar\partial v)^n.
 \end{gathered}
\end{equation*}
where the wedge product is understood in the Bedford--Taylor sense. Ko{\l}odziej and Nguyen~\cite{KN2019} extended this characterization to compact Hermitian manifolds.

Let $\Gamma\subset\C^n$ be a lattice of full real rank, and let $\mu$ be a nonzero nonnegative locally finite $\Gamma$-periodic Borel measure.
On the complex torus $\C^n/\Gamma$, the existence results in~\cite{Yau1978,Kol1998,Kol2008,DDGHKZ,KN2019} give, under~\eqref{eq:holder-measure}, a H\"older continuous $\Gamma$-periodic function $\phi_A$ satisfying $\MA(Q+\phi_A)=\mu$ for every real quadratic polynomial $Q$ with $A=\Hess Q>0$ and the mass compatibility condition
\[
 \det\Hess Q=\bar\mu:=
 \frac{\mu(\mathcal F)}{|\mathcal F|},
\]
where $\mathcal F$ is a half-open fundamental cell of $\Gamma$.

We address the corresponding classification problem: must every continuous entire psh solution of~\eqref{eq:normalization} satisfying~\eqref{eq:main-hyp} arise from a compatible quadratic background and a periodic corrector? Theorem~\ref{thm:main} gives an affirmative answer under~\eqref{eq:holder-measure}.

\begin{theorem}\label{thm:main}
Let $u\in\PSH(\C^n)\cap C^0(\C^n)$, and assume that
its complex Monge--Amp\`ere measure $\mu:=\MA(u)$ is
$\Gamma$-periodic and satisfies~\eqref{eq:holder-measure}. If
\begin{equation}\label{eq:main-hyp}
 m|z|^2-C_0\le u(z)\le M|z|^2+C_0\qquad\text{in }\C^n
\end{equation}
for some constants $m,M>0$ and $C_0\ge0$, then there is a unique decomposition
\begin{equation}\label{eq:decomposition}
 u=Q+\ell+\phi,\qquad \int_{\mathcal F}\phi\,dV=0,
\end{equation}
where $Q$ is a homogeneous real quadratic polynomial, $\ell$ is real affine, and $\phi$ is H\"older continuous and $\Gamma$-periodic. Moreover,
\begin{equation}\label{eq:matrix-main}
 m|z|^2\le Q(z)\le M|z|^2,\qquad
 A:=\Hess Q>0,\qquad \det A=\bar\mu.
\end{equation}
\end{theorem}


The measure $\mu$ need not be absolutely continuous; if it has a
density, no positive lower bound is required.


\begin{corollary}
\label{cor:holder-solution}
Let $u$ be a locally H\"older continuous plurisubharmonic function
on $\C^n$ satisfying~\eqref{eq:main-hyp}.
If $\MA(u)$ is periodic, then
\eqref{eq:decomposition}--\eqref{eq:matrix-main} hold.
\end{corollary}

\begin{proof}
Since $u$ itself provides the local subsolutions
in~\eqref{eq:holder-measure}, the conclusion follows from
Theorem~\ref{thm:main}.
\end{proof}

\begin{corollary}
\label{cor:measure-classes}
Theorem~\ref{thm:main} applies to periodic measures $\mu=f\,dV$
with $f\ge0$ and $f\in L^p(\mathcal F)$ for some $p>1$.
It also applies to periodic measures satisfying
\begin{equation}\label{eq:riesz-class}
 \mu(B_r(x))\le Kr^{2n-2+\delta}
 \qquad(0<r\le1)
\end{equation}
for some constants $K,\delta>0$.
\end{corollary}

\begin{proof}
A H\"older continuous periodic corrector exists
by~\cite{DDGHKZ} in the $L^p$ case and
by~\cite[Theorem~A]{Hiep} under~\eqref{eq:riesz-class}.
In either case, the corrector provides the local subsolutions
in~\eqref{eq:holder-measure}, so Theorem~\ref{thm:main} applies.
\end{proof}

The measure class covered by Theorem~\ref{thm:main} is larger
than the class defined by~\eqref{eq:riesz-class}.
Section~\ref{sec:singular-example} illustrates this with a purely
singular example supported on a union of real $n$-dimensional planes.
In complex dimension one, $\MA(u)=\tfrac14\Delta u$,  and the assumption~\eqref{eq:holder-measure}
can be removed; see Theorem~\ref{thm:dimension-one}.

\subsection{Outline of the proof}

To study the behavior of $u$ at infinity, we normalize the mean mass to one as in Section~\ref{subsec:periodic-conventions} and consider the blow-down family for $R\ge1$
\begin{equation}\label{eq-foundation-rescale}
 u_R(z)=R^{-2}u(Rz),\qquad
 \mu_R(E)=R^{-2n}\mu(RE),\qquad
 \MA(u_R)=\mu_R.
\end{equation}
The two-sided quadratic growth condition gives uniform local bounds for $u_R$. For each fixed $r\ge1$, let $w_R$ be the continuous psh solution of the effective problem
\begin{equation}\label{eq-foundation-replacement}
 \MA(w_R)=dV\quad\text{in }B_r,
 \qquad w_R=u_R\quad\text{on }\partial B_r.
\end{equation}
In Section~\ref{sec:homogenization}, we use a finite-ball estimate for constant-density solutions, proved in Appendix~\ref{app:estimate}, to control their deviation from local averages without involving derivatives of the boundary data. Periodic correctors yield the quantitative homogenization estimates
\[
 \norm{u_R-w_R}_{W^{1,1}(B_r)}\le CR^{-1},
 \qquad
 \norm{u_R-w_R}_{L^\infty(B_r)}\le CR^{-\gamma}
\]
for some $\gamma>0$. Together, these estimates make the blow-down family precompact in the locally uniform topology. The constant-density rigidity theorem of Hu and Sheng~\cite{HuSheng2026} then shows that every limit is a homogeneous real quadratic polynomial satisfying the averaged equation. Combining small perturbation regularity with an iteration argument (cf. Caffarelli--Li~\cite[Section~3]{CL2003}), we identify a unique quadratic part $Q$. We then derive mixed energy estimates after subtracting affine functions chosen at each scale.

In Section~\ref{sec:linear-term}, we linearize the equation around $Q$. Solving a constant-coefficient elliptic equation on the torus gives a periodic potential for $\mu-dV$. After subtracting this potential from $u-Q$, we use the mixed energy estimates to control the remaining nonlinear term and compare averaged gradients across scales. This determines the linear coefficient $b$.
The uniform homogenization estimate, combined with the coefficient estimates, also gives
\[
 \sup_{B_R}|u-Q-b\cdot x|\le CR^{2-\gamma},
 \qquad R\ge2.
\]
Singular integral estimates then yield the averaged logarithmic bound
\[
 \frac{1}{|B_R|}\int_{B_R}|u-Q-b\cdot x|\,(dV+d\mu)
 \le C\log(2+R).
\]

Finally, in Section~\ref{sec:lowpoint}, we compare $u$ with the solution $\psi=Q+b\cdot x+\phi_A$, where $A=\Hess Q$. These two functions have the same Monge--Amp\`ere measure $\mu$. Adapting the critical-density argument of Jin, Li, Tran, and Tu~\cite{JLT}, we establish a low-point estimate of weak Harnack type at large scales for their difference. 
Combining the large-scale low-point estimate with the logarithmic bound in $\mu$-average, we obtain a pointwise logarithmic bound on $u-\psi$ and then an oscillation contraction.
Together, these force the difference to be constant, completing the decomposition.

\section{Preliminaries and the one-dimensional case}
\label{subsec:periodic-conventions}\label{sec:singular-example}
 
We first note that the positive quadratic lower bound in~\eqref{eq:main-hyp} implies $\bar\mu>0$. Indeed, suppose that $\bar\mu=0$. Then $\mu\equiv 0$. Let $m>0$ and $C_0\ge0$ be the constants in~\eqref{eq:main-hyp}, and for $R>0$ write $B_R=\{z\in\C^n:|z|<R\}$. On $B_R$ set
\[
 q_R(z)=m|z|^2/2+mR^2/2-C_0.
\]
The lower growth bound and comparison give $u(0)\ge q_R(0)=mR^2/2-C_0$ for every $R$, a contradiction.
Therefore, in the subsequent sections we replace $u$ and $\mu$ by
$\bar\mu^{-1/n}u$ and $\bar\mu^{-1}\mu$, respectively, and assume that
\[
\bar\mu=1.
\]
This normalization will be used without further mention, and we retain
the same notation for the rescaled quantities. 

Next, we recall the existence of periodic correctors. Let $Q$ be a real quadratic polynomial with $A=\Hess Q>0$. We seek a $\Gamma$-periodic function $\phi_A$ such that $Q+\phi_A$ is plurisubharmonic and
\begin{equation}\label{eq:periodic-corrector-general}
 \MA(Q+\phi_A)=\mu,\qquad
 \int_{\mathcal F}\phi_A\,dV=0,
\end{equation}
where the zero-mean normalization fixes the additive constant.
For smooth $\phi_A$, periodicity and Stokes' theorem give, for $1\le k\le n$,
\[
 \int_{\C^n/\Gamma}(\dd\phi_A)^k\wedge(\dd Q)^{n-k}
 =\int_{\C^n/\Gamma}\dd\bigl(\phi_A(\dd\phi_A)^{k-1}\wedge(\dd Q)^{n-k}\bigr)=0.
\]
Consequently,
\[
 \mu(\mathcal F)
 =\frac{1}{2^n n!}\int_{\C^n/\Gamma}(\dd Q+\dd\phi_A)^n
 =\frac{1}{2^n n!}\int_{\C^n/\Gamma}(\dd Q)^n
 =(\det A)|\mathcal F|.
\]
Thus $\det A=\bar\mu$ is the mass compatibility condition. The same identity holds for continuous psh potentials by periodic smoothing and Bedford--Taylor continuity on the torus. More precisely, applying~\cite[Theorem~4.3]{DDGHKZ} to $\C^n/\Gamma$ gives the following corrector lemma.

\begin{lemma}\label{lem:measure-corrector}
Assume that $\mu$ satisfies~\eqref{eq:holder-measure}.
Let $Q$ be a real quadratic polynomial with $A=\Hess Q>0$ and
$\det A=\bar\mu$. There is a unique H\"older continuous
$\Gamma$-periodic function $\phi_A$ such that $Q+\phi_A$ is
plurisubharmonic and satisfies~\eqref{eq:periodic-corrector-general}.
The function $\phi_A$ depends only on $A$.
\end{lemma}

Condition~\eqref{eq:holder-measure} also allows purely singular periodic measures. For example, the locally Lipschitz psh function
\[
 u(z)=\sum_{j=1}^n\bigl(q(x_j)+q(y_j)\bigr),\qquad z_j=x_j+iy_j,
\]
where $q(t)=(2k+1)t-k(k+1)$ for $k\le t\le k+1$, $k\in\mathbb Z$, satisfies
\[
 \MA(u)=\mu=\frac1{2^n}\bigotimes_{j=1}^n\left(
 \sum_{k\in\mathbb Z}\delta_k(dx_j)\,dy_j
 +dx_j\sum_{l\in\mathbb Z}\delta_l(dy_j)\right).
\]
Here $q(t)-t^2$ is bounded and $1$-periodic. Moreover, $\mu(B_r(0))=\omega_n r^n$ for $0<r<1$, where $\omega_n$ is the volume of the unit ball in $\R^n$. Thus this measure fails~\eqref{eq:riesz-class} when $n\ge2$.

In complex dimension one, $\MA(u)=\tfrac14\Delta u$, and the decomposition follows from the periodic Poisson equation and the Liouville theorem for harmonic functions, without assumption~\eqref{eq:holder-measure}. The periodic Green kernel gives a mean-zero potential $\phi\in L^1(\C/\Gamma)$ satisfying $\Delta\phi=4(\mu-\bar\mu\,dV)$. Subtracting $\bar\mu|z|^2+\phi$ from $u$ gives an entire harmonic function with at most quadratic growth in mean, hence a polynomial of degree at most two. Thus $\phi$ has a continuous periodic representative, since it is the difference between $u$ and a polynomial. The bounds and mass identity for the quadratic part follow by scaling and taking the Laplacian; uniqueness follows because a periodic polynomial is constant and the periodic part has zero mean.

\begin{theorem}\label{thm:dimension-one}
Let $u\in\PSH(\C)\cap C^0(\C)$, and assume that
its complex Monge--Amp\`ere measure $\mu:=\MA(u)$ is
$\Gamma$-periodic for a lattice $\Gamma\subset\C$ of full real rank.
If~\eqref{eq:main-hyp} holds for some constants $m,M>0$ and $C_0\ge0$,
then there is a unique decomposition~\eqref{eq:decomposition},
where $Q$ is a homogeneous real quadratic polynomial, $\ell$ is real
affine, and $\phi$ is continuous and $\Gamma$-periodic.
Moreover,~\eqref{eq:matrix-main} holds.
\end{theorem}
 
\section{Basic estimates in pluripotential theory}
\label{subsec:local-approximation}

Throughout the rest of the paper, we assume $n\ge2$. 
In this section, we summarize some classical identities and standard estimates from pluripotential theory; see~\cite{BT76,BT,KNStability2019}. For completeness, we present the local arguments in the smooth case; the general case of continuous plurisubharmonic potentials follows by smooth approximation.

Set $d=2n$ and identify $\C^n$ with $\R^d$. We use
$d^c=\frac{i}2(\bar\partial-\partial)$, so that $dd^c=\dd$, and write
\[
 \omega_v=dd^cv,\qquad \omega_0=dd^c|z|^2,\qquad
 c_n=2^n n!,\qquad \omega_0^n=c_n\,dV.
\]
Mixed products of locally bounded psh potentials are understood in the
Bedford--Taylor sense~\cite{BT76,BT}. We use their continuity under locally
uniform convergence and the mixed Chern--Levine--Nirenberg estimates on
compactly contained subdomains.

Let $\eta\in C_c^\infty(B_1)$ be nonnegative with $\int\eta\,dV=1$,
and set $\eta_\varepsilon(z)=\varepsilon^{-2n}\eta(z/\varepsilon)$.
For $v\in\PSH(\Omega)\cap C(\Omega)$ and $\varepsilon_k\to0$, we consider the smooth approximation
\[
 v_k=v*\eta_{\varepsilon_k}\quad\text{on }\Omega_{\varepsilon_k},
\]
where $\Omega_\varepsilon=\{z\in\Omega:\operatorname{dist}(z,\partial\Omega)>\varepsilon\}$.
Then $v_k$ are smooth psh functions and converge locally uniformly to $v$.
For a constant-density solution
$w\in\PSH(B)\cap C(\overline B)$ on a ball $B$, we instead approximate
its boundary values uniformly by $f_k\in C^\infty(\partial B)$ and solve
\[
 \MA(w_k)=dV\quad\text{in }B,\qquad
 w_k=f_k\quad\text{on }\partial B.
\]
The solutions are smooth by~\cite{CKNS}, and comparison gives
\begin{equation}\label{eq-replacement-approx}
 \norm{w_k-w}_{L^\infty(B)}
 \le\norm{f_k-w}_{L^\infty(\partial B)}\longrightarrow0.
\end{equation}
When $w=v$ on $\partial B$ and $v$ is continuous and psh on a larger
domain, we take $f_k=v_k|_{\partial B}$, where $v_k$ are local
regularizations of $v$.

Let $v,w,p_1,\ldots,p_m$ be continuous psh functions on a domain
$\Omega$, where $0\le m\le n-1$, and set
$a=v-w$ and $T=\omega_{p_1}\wedge\cdots\wedge\omega_{p_m}$.
We define
\begin{equation}\label{eq:gradient-definition}
 da\wedge d^ca\wedge T
 =\tfrac12dd^c(a^2T)-a\,dd^ca\wedge T.
\end{equation}
An empty product is understood as $1$.
For smooth potentials, this is the product rule. Under locally uniform psh approximation, Bedford--Taylor continuity and the local mass bounds for mixed currents allow us to pass to the limit in both terms on the right:
the continuous factors converge uniformly on compact sets, and
the coefficient measures have uniformly bounded total variation.
Thus \eqref{eq:gradient-definition} agrees with the smooth expression,
defines a positive current, and is weakly continuous under locally
uniform psh approximation of all the potentials.

For smooth potentials and $n\ge2$, one has
\[
 dd^c(da\wedge d^ca\wedge T)=-(dd^ca)^2\wedge T
\]
for every closed mixed current $T$ of bidegree $(m,m)$, with $0\le m\le n-2$.
Together with the identity
\[
 (\omega_v-\omega_w)^2\wedge
 \sum_{j=0}^{n-2}(n-1-j)\omega_v^j\wedge\omega_w^{n-2-j}
 =\omega_v^n-\omega_w^n
  -n(\omega_v-\omega_w)\wedge\omega_w^{n-1},
\]
this implies
\begin{equation}\label{eq:general-current-identity}
 n\,dd^c(v-w)\wedge\omega_w^{n-1}
 =\omega_v^n-\omega_w^n+dd^cS_{v,w}(v-w),
\end{equation}
where the positive $(n-1,n-1)$-form $S_{v,w}h$ is defined by
 \begin{equation}\label{eq:general-error-current}
  S_{v,w}h
  =\sum_{j=0}^{n-2}(n-1-j)\,
  dh\wedge d^ch
  \wedge\omega_v^j\wedge\omega_w^{n-2-j}.
 \end{equation}
For $n=1$, the sum in \eqref{eq:general-error-current} is zero
and the identity is immediate.
The identity extends to continuous psh potentials by approximation.

For the gradient estimate below, we follow the proof of
Ko{\l}odziej and Nguyen~\cite[Lemma~3.6]{KNStability2019}.

\begin{lemma}\label{lem-energy}
Let $r\ge1$, and let $v$ be continuous and psh on a neighborhood of
$\overline{B_{r+2}}$, with $\norm{v}_{L^\infty(B_{r+2})}\le M$.
Suppose that $w\in\PSH(B_r)\cap C(\overline{B_r})$ satisfies
$\MA(w)=dV$ in $B_r$ and $w=v$ on $\partial B_r$. Set $h=v-w$ and
\[
 I=\int_{B_r}(v-w)\,[\MA(w)-\MA(v)].
\]
Then $I\ge0$, $h\in W^{1,1}_0(B_r)$, and
\begin{equation}\label{eq-energy-coercive}
 \left(\int_{B_r}|\nabla h|\,dV\right)^2\le CI.
\end{equation}
Moreover,
\begin{equation}\label{eq-full-energy}
 \sum_{j=0}^{n-1}\int_{B_r}i\partial h\wedge\bar\partial h
 \wedge\omega_v^j\wedge\omega_w^{n-1-j}\le CI.
\end{equation}
Here $C$ depends only on $n,r,M$.
\end{lemma}

\begin{proof}
Choose smooth psh regularizations $v_k$ converging uniformly to $v$
on $\overline{B_{r+2}}$, and let $w_k$ solve the constant-density
Dirichlet problem on $B_r$ with boundary value $v_k$.
The approximation above,
in particular~\eqref{eq-replacement-approx}, gives
$w_k\to w$ uniformly on $\overline{B_r}$.

Since $\Hess v_k\ge0$ and $\det\Hess w_k=1$, the comparison principle
implies
\[
 v_k+|z|^2-r^2\le w_k\le\sup_{\partial B_r}v_k
 \quad\text{in }B_r.
\]
Thus the functions
\begin{equation}\label{eq-psh-gluings}
 A_k=v_k+\max\{|z|^2-r^2,0\},\qquad
 W_k=\begin{cases}
 w_k&\text{in }B_r,\\
 v_k+|z|^2-r^2&\text{off }B_r
 \end{cases}
\end{equation}
are psh on $B_{r+2}$ and uniformly bounded in terms of $r$ and $M$.
These gluings will also allow us to apply the localized
integration-by-parts estimate in the homogenization argument.
Choose $\chi\in C_c^\infty(B_{r+2})$ with $0\le\chi\le1$ and
$\chi=1$ near $\overline{B_r}$. Integration by parts gives
\[
 \int_{B_r}\tr\Hess w_k\,dV\le\frac14\int_{B_{r+2}}\chi\,d(\Delta W_k)\le C.
\]

For $h_k=v_k-w_k$, integration by parts and $h_k=0$ on $\partial B_r$
give
\[
 \begin{split}
 I_k:=\int_{B_r}(v_k-w_k)\,[\MA(w_k)-\MA(v_k)]
 &=\frac1{2^n n!}\sum_{j=0}^{n-1}\int_{B_r}
 i\partial h_k\wedge\bar\partial h_k
 \wedge\omega_{v_k}^j\wedge\omega_{w_k}^{n-1-j}\\
 &\ge\frac1{2^n n!}\int_{B_r}
 i\partial h_k\wedge\bar\partial h_k\wedge\omega_{w_k}^{n-1}\\
 &\ge\frac1{4n}\int_{B_r}
 \frac{|\nabla h_k|^2}{\tr\Hess w_k}\,dV,
 \end{split}
\]
where the last inequality uses
\[
 \cof(\Hess w_k)=(\Hess w_k)^{-1}
 \ge\frac{\mathrm{Id}}{\tr\Hess w_k}.
\]
Cauchy--Schwarz and the trace bound yield
\[
 \left(\int_{B_r}|\nabla h_k|\,dV\right)^2
 \le\left(\int_{B_r}\frac{|\nabla h_k|^2}{\tr\Hess w_k}\,dV\right)
       \left(\int_{B_r}\tr\Hess w_k\,dV\right)
 \le CI_k.
\] 

The gluings in~\eqref{eq-psh-gluings} converge uniformly to bounded
continuous psh functions $A,W$ that agree off $B_r$.
Their difference $H=A-W$ is the zero extension of $h$, with
$\operatorname{supp}H\subset\overline{B_r}\Subset B_{r+2}$.
Locally bounded psh functions belong to $W^{1,2}_{\mathrm{loc}}$,
so $H\in W^{1,2}_{\mathrm{loc}}(B_{r+2})$.
Its restriction $h$ has zero Sobolev trace and belongs to
$W^{1,1}_0(B_r)$.

Uniform convergence and Bedford--Taylor continuity give
\[
 I_k\longrightarrow\int_{B_r}(v-w)\,[\MA(w)-\MA(v)]=I.
\]
To justify the integral over the whole ball, the
Chern--Levine--Nirenberg estimate on a larger ball bounds the masses
of $\MA(v_k)$ and $\MA(v)$ on $B_r$ uniformly. For
$E_\delta=B_r\setminus B_{r-\delta}$, $0<\delta<r$, we have
\[
 \int_{E_\delta}|h_k|\,[\MA(v_k)+\MA(w_k)]
 \le C\left(\norm{h_k-h}_{L^\infty(B_r)}
             +\sup_{E_\delta}|h|\right).
\]
The right-hand side tends to zero as first $k\to\infty$ and then
$\delta\to0$, since $h=0$ on $\partial B_r$.
The same boundary-layer bound holds for $|h|[\MA(v)+\MA(w)]$.
Bedford--Taylor continuity with interior cutoffs therefore proves the
displayed convergence, even if $\MA(v)$ charges $\partial B_r$.
Hence $I\ge0$, and lower semicontinuity proves~\eqref{eq-energy-coercive}.

By the weak continuity of mixed gradient currents, the positive
mixed energy measures in the identity for $I_k$ converge locally
in $B_r$. Exhausting $B_r$ and using $I_k\to I$
proves~\eqref{eq-full-energy}.
\end{proof}

The following integration-by-parts estimate will be applied to the rescaled periodic potential and used to pass from Lebesgue integral bounds to bounds with respect to $\mu$.

\begin{lemma}\label{lem:local-transfer}
Let $0<r<s$ and let $v,w,a,b\in\PSH(B_s)\cap C(B_s)$ satisfy
\[
 \norm{v}_{L^\infty(B_s)}+\norm{w}_{L^\infty(B_s)}
 +\norm{a}_{L^\infty(B_s)}+\norm{b}_{L^\infty(B_s)}\le M.
\]
For a real $\chi\in C_c^\infty(B_r)$, one has
\begin{equation}\label{eq:local-transfer}
 \left|\int\chi\,(v-w)_+\,[\MA(a)-\MA(b)]\right|
 \le C\norm{a-b}_{L^\infty(B_s)},
\end{equation}
where $C$ depends only on $n,r,s,M$ and $\norm{\chi}_{C^2}$.
\end{lemma}

\begin{proof}
Put $f=(v-w)_+$ and
\[
 T=\sum_{j=0}^{n-1}\omega_a^j\wedge\omega_b^{n-1-j}.
\]
The current $T$ is closed and positive, and
\[
 \MA(a)-\MA(b)=c_n^{-1}dd^c(a-b)\wedge T.
\]
Since $f=\max\{v,w\}-w$, it is the difference of two uniformly bounded
psh functions.
We first smooth the psh functions in this representation, as well as
$a,b$, on a slightly smaller ball. The estimates below are uniform in
these regularizations, whose subscripts we suppress.

The mixed Chern--Levine--Nirenberg estimates give
\[
 \int_{B_r}\omega_0\wedge T
 +\bigl\|dd^cf\wedge T\bigr\|_{\mathrm{TV}(B_r)}\le C,
\]
where $\|\nu\|_{\mathrm{TV}(B_r)}$ denotes the total variation of the signed measure $\nu$ on $B_r$.
Choose a nonnegative cutoff $\eta\in C_c^\infty(B_s)$ equal to one
on $B_r$. The same bounds on a fixed neighborhood of its support and
\eqref{eq:gradient-definition} imply
\[
 \begin{aligned}
 \int_{B_r}df\wedge d^cf\wedge T
 &\le\int\eta\,df\wedge d^cf\wedge T\\
 &=\frac12\int f^2dd^c\eta\wedge T
   -\int\eta f\,dd^cf\wedge T\le C.
 \end{aligned}
\]
Expanding $dd^c(\chi f)\wedge T$, we use the preceding estimates to bound the total variation of the terms involving $dd^c\chi$ and $dd^cf$ by $C$. The Cauchy--Schwarz inequality bounds the total variation of each cross term by
\[
 C\left(\int d\chi\wedge d^c\chi\wedge T\right)^{1/2}
  \left(\int_{B_r}df\wedge d^cf\wedge T\right)^{1/2}\le C.
\]
Consequently,
\begin{equation}\label{eq:local-transfer-TV}
 \bigl\|dd^c(\chi f)\wedge T\bigr\|_{\mathrm{TV}}\le C.
\end{equation}
Since $\chi f$ has compact support, integration by parts gives
\[
 \int\chi f\,[\MA(a)-\MA(b)]
 =c_n^{-1}\int(a-b)\,dd^c(\chi f)\wedge T.
\]
The bound~\eqref{eq:local-transfer-TV} proves~\eqref{eq:local-transfer}
for the regularizations. Local uniform convergence, Bedford--Taylor continuity, and the uniform total variation bound allow us to pass to the limit in both the identity and the estimate. This also defines the signed measure in
\eqref{eq:local-transfer-TV} for the original continuous potentials.
\end{proof}

\section{Homogenization and blow-down limits}\label{sec:homogenization}

In this section, we prove the homogenization estimate and classify the possible blow-down limits.

Fix $r\ge1$. For each $R\ge1$, let $u_R$
denote the rescaled solution defined in~\eqref{eq-foundation-rescale}.
The growth condition~\eqref{eq:main-hyp} gives
\begin{equation}\label{eq-rescaled-growth}
 m|z|^2-C_0R^{-2}\le u_R\le M|z|^2+C_0R^{-2}.
\end{equation}
Let $w_R\in\PSH(B_r)\cap C(\overline{B_r})$ be the unique
solution of the effective problem~\eqref{eq-foundation-replacement},
whose existence and uniqueness follow from
\cite[Theorem~D]{BT76}.
Let $h_R=u_R-w_R$. Lemma~\ref{lem-energy}, applied to
$(v,w)=(u_R,w_R)$, gives $h_R\in W^{1,1}_0(B_r)$ and the energy estimates
\eqref{eq-energy-coercive}--\eqref{eq-full-energy},
with constants uniform in $R$ by~\eqref{eq-rescaled-growth}.

To prove the homogenization estimates, let $Q=|z|^2$ and let $\phi=\phi_{\Id}$ be the periodic corrector
given by Lemma~\ref{lem:measure-corrector}. For $R\ge1$, define
\begin{equation}\label{eq:measure-rescaled-potential}
 \phi_R(z)=R^{-2}\phi(Rz),\qquad \psi_R=Q+\phi_R,
 \qquad \MA(\psi_R)=\mu_R,
 \qquad \norm{\phi_R}_\infty\le CR^{-2}.
\end{equation}
Set
\[
 T_R=\sum_{j=0}^{n-1}\omega_{\psi_R}^j\wedge\omega_{Q}^{n-1-j},
 \qquad
 \mu_R-dV=\frac{\dd\phi_R\wedge T_R}{2^n n!}.
\]
The current $T_R$ is closed and positive.

\begin{proposition}\label{prop-homogenization}
On each fixed ball $B_r$ with $r\ge1$, we have
\begin{equation}\label{eq-homogenization-energy}
 -\int_{B_r}(u_R-w_R)\,(d\mu_R-dV)\le C_rR^{-2}.
\end{equation}
Moreover,
\begin{equation}\label{eq-homogenization-W11}
 \norm{u_R-w_R}_{W^{1,1}(B_r)}\le C_rR^{-1},
\end{equation}
and we have
\begin{equation}\label{eq:measure-weighted-error}
 \int_{B_r}|u_R-w_R|\,d\mu_R\le C_rR^{-1}.
\end{equation}
\end{proposition}
\begin{proof}
Let $A,W$ be the potentials constructed in the proof of
Lemma~\ref{lem-energy} for $(v,w)=(u_R,w_R)$.
By~\eqref{eq-rescaled-growth}, they are uniformly bounded on
$B_{r+2}$.
Choose a cutoff $\chi$ supported in $B_{r+1}$ and equal to one
near $\overline{B_r}$.
The potentials $\psi_R,Q$ are also uniformly bounded on $B_{r+2}$.

Extend $h_R$ by zero outside $B_r$, so that $h_R=A-W$
on $B_{r+2}$.
Apply Lemma~\ref{lem:local-transfer} to $(A,W,\psi_R,Q)$ and
$(W,A,\psi_R,Q)$. Using $h_R=(h_R)_+-(-h_R)_+$, we obtain
\[
 \begin{aligned}
 0&\le -\int_{B_r}h_R\,(d\mu_R-dV)\\
 &=-\int\chi h_R\,[\MA(\psi_R)-\MA(Q)]
 \le C_r\norm{\phi_R}_\infty\le C_rR^{-2}.
 \end{aligned}
\]
This proves~\eqref{eq-homogenization-energy}. By~\eqref{eq-energy-coercive},
\[
 \int_{B_r}|\nabla h_R|\,dV
 \le C_r\left[-\int_{B_r}h_R\,(d\mu_R-dV)\right]^{1/2}
 \le C_rR^{-1}.
\]
Since $h_R\in W^{1,1}_0(B_r)$, Poincar\'e's inequality gives
\eqref{eq-homogenization-W11}.

Adding the two estimates and using $|h_R|=(h_R)_++(-h_R)_+$ gives
\[
 \left|\int_{B_r}|h_R|\,(d\mu_R-dV)\right|
 \le C_r\norm{\phi_R}_\infty\le C_rR^{-2}.
\]
Combining this with~\eqref{eq-homogenization-W11} gives
\[
 \begin{aligned}
 \int_{B_r}|h_R|\,d\mu_R
 &\le \int_{B_r}|h_R|\,dV
 +\left|\int_{B_r}|h_R|\,(d\mu_R-dV)\right|\\
 &\le C_rR^{-1}+C_rR^{-2}\le C_rR^{-1},
 \end{aligned}
\]
which proves~\eqref{eq:measure-weighted-error}.
\end{proof}

\begin{lemma}\label{lem:uniform-homogenization}
For each fixed $r\ge1$,
\begin{equation}\label{eq-homogenization-uniform}
 \norm{u_R-w_R}_{L^\infty(B_r)}
 \le C_rR^{-1/(2n+1)}
 \qquad(R\ge1).
\end{equation}
\end{lemma}
\begin{proof}
We combine the integral estimates in
Proposition~\ref{prop-homogenization} with a uniform capacity bound.
Recall the reference potentials
$\psi_R=Q+\phi_R$ from~\eqref{eq:measure-rescaled-potential}.
If $\phi$ has H\"older exponent $\alpha\in(0,1]$, then its rescaling satisfies
\[
 [\phi_R]_{C^{0,\alpha}}\le CR^{\alpha-2}.
\]
Hence $\psi_R=Q+\phi_R$ is uniformly bounded in $C^{0,\alpha}$ on every fixed ball.

For a Borel set $E\subset B_r$, we use the relative Bedford--Taylor capacity
\[
 \operatorname{Cap}_{B_r}(E)=\sup\left\{
 \int_E(i\partial\bar\partial\eta)^n:
 \eta\in\PSH(B_r),\ -1\le\eta\le0\right\}.
\]
The capacity is defined with respect to the unnormalized
Monge--Amp\`ere measure; consequently, the test-function integrals
in its definition differ by the factor $2^n n!$ from those computed
using $\MA$. The capacity estimate for Monge--Amp\`ere measures with H\"older continuous potentials
\cite[Proposition~3.1]{NguyenII} gives, for every $\tau>0$,
\begin{equation}\label{eq:measure-capacity}
 \mu_R(E)\le C_{r,\tau}\operatorname{Cap}_{B_r}(E)^{1+\tau}.
\end{equation}
To reduce to the zero-boundary setting, we work on a fixed larger
ball and replace $\psi_R$ by $\max\{\psi_R-C,B\rho\}$,
where $\rho$ is a strictly psh defining function of the larger ball.
For suitable uniform constants $C$ and $B$, this modified potential
agrees with $\psi_R-C$ near $\overline{B_r}$ and with $B\rho$
near the outer boundary.
Since relative capacity decreases as the domain increases,
this proves~\eqref{eq:measure-capacity} for every Borel set
$E\subset B_r$.

We first establish a one-sided stability estimate that will be applied to $u_R$ and $w_R$. Let $a,b\in\PSH(B_r)\cap C(\overline{B_r})$ satisfy $a\ge b$ on $\partial B_r$, and write $\nu=\MA(a)$. Assume that
\[
 \nu(E)\le C\operatorname{Cap}_{B_r}(E)^{1+\tau}.
\]
Set $J=\int_{B_r}(b-a)_+\,d\nu$.
For $\eta\in\PSH(B_r)$ with $-1\le\eta\le0$, we have
\[
 \{a<b-s-t\}\subset\{a<b+t\eta-s\}\subset\{a<b-s\}.
\]
Using the comparison principle and
$\MA(b+t\eta-s)\ge t^n\MA(\eta)$, and then taking the supremum
over admissible $\eta$, we obtain
\begin{equation}\label{eq:measure-stability-cap}
 t^n\operatorname{Cap}_{B_r}\{a<b-s-t\}
 \le C_n\nu\{a<b-s\}\le C_nJ/s
 \qquad(s,t>0).
\end{equation}

Write $F(s)=\operatorname{Cap}_{B_r}\{a<b-s\}^{1/n}$.
Combining the first inequality in~\eqref{eq:measure-stability-cap}
with the capacity bound for $\nu$ and taking $n$-th roots gives
\[
 tF(s+t)\le C F(s)^{1+\tau}.
\]
Starting from $s_0=s$, set $s_{j+1}=s_j+2C F(s_j)^\tau$.
Then $F(s_j)\le2^{-j}F(s)$, and summing the geometric series gives
$\sum_{j\ge0}(s_{j+1}-s_j)\le C F(s)^\tau$.
Since $F$ is nonincreasing, $F(\lim_{j\to\infty}s_j)=0$.
We therefore obtain
\[
 \sup_{B_r}(b-a)
 \le s+C\operatorname{Cap}_{B_r}\{a<b-s\}^{\tau/n}.
\]
Taking $s=t=\varepsilon$ in~\eqref{eq:measure-stability-cap}
and substituting the resulting capacity bound into the preceding
estimate at $2\varepsilon$ gives
\begin{equation}\label{eq:measure-stability-power}
 \sup_{B_r}(b-a)
 \le 2\varepsilon+C(J\varepsilon^{-n-1})^{\tau/n}
 \le C J^{\tau/[n+(n+1)\tau]},
\end{equation}
where the last inequality uses
$\varepsilon=J^{\tau/[n+(n+1)\tau]}$ for $J>0$.
The case $J=0$ follows by letting $\varepsilon\to0$.

We now apply~\eqref{eq:measure-stability-power} with
$a=u_R$, $b=w_R$, $\nu=\mu_R$, and $\tau=1$.
The capacity bound is~\eqref{eq:measure-capacity}, while
\eqref{eq:measure-weighted-error} gives $J\le C_rR^{-1}$.
Hence $\sup_{B_r}(w_R-u_R)\le C_rR^{-1/(2n+1)}$.
For the reverse order, take $a=w_R$, $b=u_R$, and $\nu=dV$.
Lebesgue measure satisfies the required capacity bound, and
\eqref{eq-homogenization-W11} gives $J\le C_rR^{-1}$, yielding
the same estimate for $\sup_{B_r}(u_R-w_R)$. Thus
\[
 \sup_{B_r}|u_R-w_R|\le C_rR^{-1/(2n+1)},
\]
which proves~\eqref{eq-homogenization-uniform}.
\end{proof}

Next, using the homogenization estimate, we identify the blow-down limit as a quadratic polynomial and prove its uniqueness.

Set $\rho(z)=1+|z|^2$. For $t>0$, define the ball average
\[
 \B_t w(z)=\frac1{|B_t|}\int_{B_t}w(z+h)\,dh,
\]
and its weighted average
\begin{equation}\label{eq:G}
 \G_t w(z)=\frac2{t^2}\int_0^t s\,\B_s w(z)\,ds.
\end{equation}
For smooth $w$, symmetry gives
\[
 \G_t w-w=\frac{t^2}{8(n+1)}\Delta_{\R}w+o(t^2).
\]
The kernel in \eqref{eq:G} is positive, has integral one, and is supported
in $\overline B_t$.

\begin{proposition}\label{prop:finite}
Let $L\ge16$ and let $w\in\PSH(B_{2L})\cap C(\overline{B_{2L}})$.
Suppose $\MA(w)=dV$ and $m_0\rho\le w\le M_0\rho$,
where $0<m_0\le M_0$. Define
\begin{equation}\label{eq:finite-constants}
 \lambda_*=\frac{m_0^{n-1}}{(4n)^{n-1}},\qquad
 t_0=\min\{1,\sqrt{m_0\lambda_*/2}\},\qquad
 \beta=\frac{m_0}{16M_0+m_0}.
\end{equation}
Then, for $|z|<\sqrt L$ and $0<t\le t_0$,
\begin{equation}\label{eq:finite}
 0\le\G_t w(z)-w(z)
 \le \frac{t^2}{\lambda_*}+\frac{4M_0^2}{m_0}L^{-\beta}\rho(z).
\end{equation}
\end{proposition}

\begin{proof}
The smooth case is proved in Appendix~\ref{app:estimate}.
For continuous $w$, the Dirichlet approximation in
Section~\ref{subsec:local-approximation} gives smooth solutions $w_k$
with $\MA(w_k)=dV$ and $w_k\to w$ uniformly on $\overline{B_{2L}}$.
Fix $0<\varepsilon<m_0$. For sufficiently large $k$,
$(m_0-\varepsilon)\rho\le w_k\le(M_0+\varepsilon)\rho$.
Apply the smooth estimate with these constants and let $k\to\infty$.
For each $0<t<t_0$, letting $\varepsilon\to0$ gives
\eqref{eq:finite} with the constants in~\eqref{eq:finite-constants}.
Continuity in $t$ gives the endpoint $t=t_0$.
\end{proof}

\begin{proposition}\label{prop-blowdown}
As $R\to\infty$, the family $u_R$ is precompact in
$C^0_{\rm loc}(\C^n)$, and for every $K>0$,
\[
 \lim_{R\to\infty}\inf_{Q\in\mathcal Q}
 \norm{u_R-Q}_{L^\infty(B_K)}=0,
\]
where $\mathcal Q$ is the compact family of homogeneous real
quadratic polynomials $Q$ satisfying
\begin{equation}\label{eq-quadratic-family}
 m|z|^2\le Q(z)\le M|z|^2,
 \qquad \det\Hess Q=1.
\end{equation}
\end{proposition}
\begin{proof}
For each fixed $L\ge16$, let $w_R$ solve the effective problem on $B_{2L}$.
By \eqref{eq-homogenization-uniform}, it differs from $u_R$ by at
most one once $R$ is sufficiently large. The bounds in
\eqref{eq-rescaled-growth} therefore give constants $C_1,m_0,M_0>0$,
independent of $L$, such that
\[
 m_0(1+|z|^2)\le w_R(z)+C_1\le M_0(1+|z|^2)
 \quad\text{on }B_{2L}.
\]
For example one can take $C_1=C_0+2$ and slightly enlarge the upper
constant. Apply Proposition~\ref{prop:finite} directly to $w_R+C_1$,
with constants independent of $L$ and $R$, and transfer the estimate
to $u_R$ using~\eqref{eq-homogenization-uniform}.
For each fixed $K$ and $t$, first take the limsup as $R\to\infty$
and then let $L\to\infty$. We obtain
\begin{equation}\label{eq-asymptotic-mean-new}
 \limsup_{R\to\infty}
 \norm{\G_t u_R-u_R}_{L^\infty(B_K)}\le Ct^2,
 \qquad 0<t\le t_0.
\end{equation}
Moreover, $\G_tu_R-u_R\ge0$ by subharmonicity.

For fixed $t$, the kernel of $\G_t$ is continuous under translations
in $L^1$. The local bounds in \eqref{eq-rescaled-growth} therefore
make $\G_tu_R$ equicontinuous on $B_K$, uniformly in $R$.
Choosing $t$ small in \eqref{eq-asymptotic-mean-new} then gives
asymptotic equicontinuity of $u_R$.
The same local bounds and a diagonal argument show that every
sequence $R_j\to\infty$ has a subsequence along which $u_{R_j}$
converges locally uniformly to a function $U$.

Along the same subsequence, the effective solutions on each fixed
ball converge uniformly to $U$ by \eqref{eq-homogenization-uniform}.
Bedford--Taylor continuity, \eqref{eq-foundation-replacement},
and \eqref{eq-rescaled-growth} yield
\begin{equation}\label{eq-U-equation-new}
 \mathcal M(U)=dV,\qquad
 m|z|^2\le U(z)\le M|z|^2.
\end{equation}
By~\eqref{eq-U-equation-new}, $\MA(U+1)=dV$ and
\[
 \min\{1,m\}(1+|z|^2)\le U(z)+1
 \le\max\{1,M\}(1+|z|^2).
\]
We now show directly that $U$ is quadratic; this also follows from Hu and Sheng's theorem~\cite[Theorem~1.1]{HuSheng2026}, applied to $U+1$.

Passing to the limit in \eqref{eq-asymptotic-mean-new} gives
$0\le\G_tU-U\le Ct^2$. For a nonnegative compactly supported
smooth test function $\eta$, symmetry of the averaging kernel gives
\[
 0\le\int U\frac{\G_t\eta-\eta}{t^2}\,dV
 =\int\eta\frac{\G_tU-U}{t^2}\,dV\le C\int\eta\,dV.
\]
The Taylor expansion of this averaging operator is
$\G_t\eta-\eta=t^2\Delta\eta/[8(n+1)]+o(t^2)$.
The preceding inequality therefore proves, distributionally,
\begin{equation}\label{eq-U-laplacian-new}
 0\le\Delta U\le C.
\end{equation}

Poisson regularity gives $U\in W^{2,p}_{\rm loc}$ for every finite
$p$. The interior regularity theorem of B{\l}ocki--Dinew \cite{BD},
applied to \eqref{eq-U-equation-new}, makes $U$ smooth; in complex
dimension one this is just Poisson regularity. Equations
\eqref{eq-U-equation-new} and \eqref{eq-U-laplacian-new} now give
\begin{equation}\label{eq-U-uniform-ellipticity-new}
 \lambda\mathrm{Id}\le\Hess U\le\Lambda\mathrm{Id}
 \quad\text{on }\C^n.
\end{equation}

To apply the real Evans--Krylov estimate, let $\mathfrak p(N)$ denote
the complex Hessian associated with a real symmetric matrix $N$. For $N\ge0$,
$\mathfrak p(N)\ge0$ and $\tr\mathfrak p(N)=\tr N/4$.
The concave extension
\[
 \mathcal E(N)=\inf_{\lambda\mathrm{Id}\le H\le\Lambda\mathrm{Id}}
 \{\log\det H+\tr[H^{-1}(\mathfrak p(N)-H)]\}
\]
is uniformly elliptic on real symmetric matrices and agrees with
$\log\det\mathfrak p(N)$ in the matrix range
\eqref{eq-U-uniform-ellipticity-new}. Thus
$\mathcal E(\HR U)=0$.
Apply the Evans--Krylov estimate \cite{Evans} to
$s^{-2}U(s\,\cdot)$ on $B_2$, using \eqref{eq-U-equation-new}:
\[
 [\HR U]_{C^\alpha(B_s)}\le Cs^{-\alpha}.
\]
Letting $s\to\infty$ makes $\HR U$ constant. The bounds in
\eqref{eq-U-equation-new} eliminate the constant and linear parts,
so $U$ satisfies \eqref{eq-quadratic-family}.

This argument applies to every subsequential limit. Since
$\mathcal Q$ is compact, precompactness proves the asymptotic
statement by contradiction.
\end{proof}

Furthermore, combining a standard iteration argument (see Caffarelli--Li~\cite{CL2003}) with interior regularity estimates, we obtain the following.

\begin{proposition}\label{prop-quadratic-approx}
There is a unique homogeneous real quadratic $Q$ satisfying
\eqref{eq-quadratic-family} such that $u_R\to Q$ locally uniformly.
For all sufficiently large $R$ there is a real affine $\ell_R$ such that
\begin{equation}\label{eq-affine-L2-approx}
 \left(\fint_{B_R}|u-Q-\ell_R|^2\,dV\right)^{1/2}\le CR,
\end{equation}
and, for $0\le j\le n-1$,
\begin{equation}\label{eq-mixed-energy}
 \int_{B_R}i\,\partial(u-Q-\ell_R)
 \wedge\bar\partial(u-Q-\ell_R)
 \wedge\omega_u^{j}\wedge\omega_Q^{n-1-j}\le CR^{2n}.
\end{equation}
\end{proposition}
\begin{proof} 
Let $w_R$ solve the effective problem on $B_4$.
By Proposition~\ref{prop-blowdown} and
\eqref{eq-homogenization-uniform}, its $L^\infty(B_4)$ distance
from the compact quadratic family~\eqref{eq-quadratic-family}
tends to zero.
The constant-density perturbation estimate of Hu and
Sheng~\cite[Lemma~6.1]{HuSheng2026}, based on Savin's
small-perturbation regularity theorem~\cite{Savin}, therefore
gives uniform $C^{2,\alpha}$ estimates and ellipticity bounds
on $B_3$ by application on fixed interior balls.
Interior Schauder estimates then yield, for all sufficiently large $R$,
\begin{equation}\label{eq-replacement-smooth-interior}
 \lambda_1\mathrm{Id}\le\Hess w_R\le\Lambda_1\mathrm{Id},
 \qquad \norm{w_R}_{C^k(B_2)}\le C_k.
\end{equation}

Let $\mathcal P$ denote the class of real quadratic polynomials $P$ with
$\Hess P>0$ and $\det\Hess P=1$, and set

\begin{equation}\label{eq-excess-definition}
 E(R)=R^{-2}\inf_{P\in\mathcal P}\fint_{B_R}|u-P|\,dV.
\end{equation}
The infimum has a minimizer $P_R$. Indeed, a minimizing sequence
has bounded rescaled coefficients by equivalence of polynomial
norms on $B_1$; its complex Hessians have bounded largest eigenvalue
and determinant one, so their smallest eigenvalues stay
positive. The constraint is therefore preserved in the limit.

Write $p_R(z)=R^{-2}P_R(Rz)$. The same norm equivalence, the size
bound \eqref{eq-rescaled-growth}, and comparison with any fixed
member of $\mathcal P$ show
\begin{equation}\label{eq-polynomial-coefficient-bounds}
 E(R)\le C,\qquad
 \norm{p_R}_{C^2(B_2)}\le C,\qquad
 \lambda_2\mathrm{Id}\le\Hess p_R\le\Lambda_2\mathrm{Id}.
\end{equation}
Since $\det\Hess w_R=\det\Hess p_R$, the difference
$w_R-p_R$ solves the homogeneous linear equation obtained by
integrating the differential of $\log\det$ along the segment
between their complex Hessians. By
\eqref{eq-replacement-smooth-interior} and
\eqref{eq-polynomial-coefficient-bounds}, its real coefficient
matrix is uniformly elliptic with uniform interior $C^k$ bounds.
Interior estimates for this linear equation consequently give
\begin{equation}\label{eq-linear-C3-approx}
 \norm{w_R-p_R}_{C^3(B_{1/2})}
 \le C\norm{w_R-p_R}_{L^1(B_1)}
 \le C\bigl(E(R)+R^{-1}\bigr),
\end{equation}
where the last inequality is \eqref{eq-homogenization-W11} and
\eqref{eq-excess-definition}.

Let $T_R$ be the quadratic Taylor polynomial of $w_R$ at the
origin. Its complex Hessian has determinant one, so the
rescaled polynomial $R^2T_R(\,\cdot/R)$ is in $\mathcal P$.
Taylor's theorem in \eqref{eq-linear-C3-approx}, together with
\eqref{eq-homogenization-W11}, gives for $0<\theta<1/4$
\begin{equation}\label{eq-excess-improvement}
 E(\theta R)\le C\theta E(R)+C\theta^{-d-2}R^{-1}.
\end{equation}
In this estimate the term $C\theta E(R)$ is the normalized cubic
Taylor remainder; restricting the $L^1$ error between $u_R$ and the effective solution from
$B_1$ to $B_\theta$ costs $\theta^{-d-2}$.

Choose $\theta$ so that $a=C\theta<1$, and denote the last fixed
constant in \eqref{eq-excess-improvement} by $B$.
For a fixed sufficiently large $r$, iterate that inequality from
radius $\theta^{-k}r$ down to radius $r$:
\[
 E(r)\le a^k E(\theta^{-k}r)
       +\frac{B\theta}{r}\sum_{j=0}^{k-1}(a\theta)^j.
\]
The uniform bound~\eqref{eq-polynomial-coefficient-bounds} gives
$a^kE(\theta^{-k}r)\to0$ as $k\to\infty$.
Summing the geometric series yields
\begin{equation}\label{eq-optimal-L1-rate}
 E(R)\le C/R,\qquad
 \fint_{B_R}|u-P_R|\,dV\le CR.
\end{equation}

In real coordinates write
\begin{equation}\label{eq-PR-coefficients}
 P_R(x)=\tfrac12x^T H_Rx+\beta_R\cdot x+c_R.
\end{equation}
Compare $P_R$ and $P_{2R}$ on $B_R$ using
\eqref{eq-optimal-L1-rate}. Polynomial norm equivalence gives
\[
 |H_{2R}-H_R|\le C/R,\qquad
 |\beta_{2R}-\beta_R|\le C,\qquad
 |c_{2R}-c_R|\le CR.
\]
The quadratic increments are summable along dyadic radii, giving a
limit $H_\infty$. The same comparison on $B_R$ for $R\le S\le2R$
yields, for all sufficiently large $R$,
\begin{equation}\label{eq-quadratic-coefficient-rate}
 |H_R-H_\infty|\le C/R.
\end{equation}
The other two inequalities give $\beta_R=O(\log R)$ and $c_R=O(R)$.

It follows from \eqref{eq-PR-coefficients} and
\eqref{eq-optimal-L1-rate} that $u_R$ converges in $L^1(B_1)$ to
$Q(x)=x^TH_\infty x/2$. Proposition~\ref{prop-blowdown} turns this
into locally uniform convergence and proves
\eqref{eq-quadratic-family}. This also proves uniqueness of $Q$.

Combining \eqref{eq-linear-C3-approx} and \eqref{eq-optimal-L1-rate} gives
\begin{equation}\label{eq-C1-optimal-rate}
 \norm{w_R-p_R}_{C^1(B_{1/2})}\le C/R.
\end{equation}
For the effective solution in the original variables, $\widehat w_R(x)=R^2w_R(x/R)$,
let $a_R(x)=\beta_R\cdot x+c_R$.
Equations~\eqref{eq-quadratic-coefficient-rate} and
\eqref{eq-C1-optimal-rate} imply
\begin{equation}\label{eq-physical-gradient}
 |\nabla(\widehat w_R-Q-a_R)|\le C
 \quad\text{on }B_{R/2}.
\end{equation}
The full energy estimate \eqref{eq-full-energy} and
\eqref{eq-homogenization-energy} scale by the factor $R^{d+2}$.
Using the uniform lower bound for $\Hess w_R$ in
\eqref{eq-replacement-smooth-interior}, we obtain, for each
$0\le j\le n-1$,
\begin{equation}
 \int_{B_{R/2}}i\,\partial(u-\widehat w_R)
 \wedge\bar\partial(u-\widehat w_R)
 \wedge\omega_u^{j}\wedge\omega_0^{n-1-j}\le CR^d.
\end{equation}
Here $\omega_0=\dd|z|^2$. Every mixed term is individually bounded because all terms in
\eqref{eq-full-energy} are positive.

The Chern--Levine--Nirenberg cutoff estimate and the growth bound
\eqref{eq:main-hyp}, rescaled to fixed balls, give
\begin{equation}\label{eq-growth-mixed-mass}
 \int_{B_R}\omega_u^{j}\wedge\omega_0^{n-j}\le CR^d
 \quad(0\le j\le n).
\end{equation}
We write
\[
 u-Q-a_R=(u-\widehat w_R)+(\widehat w_R-Q-a_R).
\]
Using the positivity of the mixed currents and the inequality $|\xi+\eta|^2\le2|\xi|^2+2|\eta|^2$, we combine~\eqref{eq-physical-gradient}--\eqref{eq-growth-mixed-mass} to obtain~\eqref{eq-mixed-energy} on $B_{R/2}$, initially with
$\omega_0$ instead of $\omega_Q$; the two fixed positive forms
are comparable by \eqref{eq-quadratic-family}.

Finally replace $R$ by $2R$ and set $\ell_R=a_{2R}$.
Equations~\eqref{eq-optimal-L1-rate} and
\eqref{eq-quadratic-coefficient-rate} give
$\fint_{B_R}|u-Q-\ell_R|\le CR$.
The $j=0$ case of \eqref{eq-mixed-energy} bounds its gradient in
$L^2(B_R)$ by $CR^{d/2}$. Poincar\'e's inequality, with the
preceding $L^1$ bound controlling the mean, proves
\eqref{eq-affine-L2-approx}.
\end{proof}

\section{The linear term in the asymptotic expansion}\label{sec:linear-term}

Proposition~\ref{prop-quadratic-approx} allows the affine correction
$\ell_R$ to vary with $R$. In this section, we first identify a fixed vector $b$ and
establish uniform averaged energy bounds
for $u-Q-b\cdot x$.

\begin{definition}
Let $A=\Hess Q$, with $\cof A=(\det A)A^{-1}$.
Define the constant-coefficient uniformly elliptic operator
\[
 L_Aa=c_n\tr\bigl((\cof A)\Hess a\bigr).
\]
Then $n\,dd^ca\wedge\omega_Q^{n-1}=(L_Aa)\,dV$.
A distribution $H$ is called $A$-harmonic if $L_AH=0$.
\end{definition}

Set $h=u-Q$. For a real affine function $\ell$, we abbreviate
\[
 S(h-\ell):=S_{u,Q}(h-\ell),
\]
where $S_{u,Q}$ is defined in~\eqref{eq:general-error-current}.
Applying~\eqref{eq:general-current-identity} to
$(v,w)=(u-\ell,Q)$ gives
\begin{equation}\label{eq:current-identity}
 L_Ah\,dV=c_n(\mu-dV)+dd^cS(h-\ell).
\end{equation}
In particular, the current $dd^cS(h-\ell)$ is independent of $\ell$.
The mixed energy estimate~\eqref{eq-mixed-energy} and the
comparability of $\omega_Q$ and $\omega_0$ imply
\begin{equation}\label{eq:local-current-mass}
 \int_{B_R}S(h-\ell_R)\wedge\omega_0\le CR^d.
\end{equation}
Here and below, passing to fixed enlargements of a ball only changes the constant $C$.

Since $\mu-dV$ has zero mass on $\C^n/\Gamma$, there exists $p\in L^1(\C^n/\Gamma)$ such that
\[
 L_Ap=c_n(\mu-dV),\qquad \int_{\C^n/\Gamma}p\,dV=0,
\]
where the equation holds in the sense of distributions. Identifying $p$ with its periodic extension, we have
\begin{equation}\label{eq:periodic-L1-potential}
 \int_{B_R}|p|\,dV\le CR^d\qquad(R\ge1).
\end{equation}
Thus~\eqref{eq:current-identity} becomes
\begin{equation}\label{eq:linear-error}
 L_A(h-p)\,dV=dd^cS(h-\ell).
\end{equation}

\begin{lemma}\label{lem:affine}
There is $b\in\R^d$ such that, with $g_0=h-b\cdot x$,
the following estimates hold for $R\ge2$:
\begin{equation}\label{eq:fixed-energy}
 \fint_{B_R}|g_0|^2\le CR^2.
\end{equation}
Moreover,
\begin{equation}\label{eq:fixed-gradient}
 \int_{B_R}|\nabla g_0|^2\le CR^d,
\end{equation}
and we have
\begin{equation}\label{eq:fixed-current}
 \int_{B_R}S(g_0)\wedge\omega_0\le CR^d.
\end{equation}
\end{lemma}
\begin{proof} 
Choose a nonnegative even probability kernel $\rho\in C_c^\infty$ that is
radial in coordinates reducing $L_A$ to the Laplacian. Put
$\rho_R(x)=R^{-d}\rho(x/R)$ and
\[
 b(R)=\nabla(\rho_R*h)(0),\qquad a(R)=(\rho_R*h)(0).
\]
To compare the averages at scales $R$ and $2R$, we solve a Poisson
equation for $\rho_{2R}-\rho_R$. Its radiality and zero integral
give a solution $\psi_R$ supported in a fixed multiple of $B_R$:
\[
 L_A\psi_R=\rho_{2R}-\rho_R,\qquad
 |D^k\psi_R|\le C_kR^{2-d-k}\quad(k=2,3).
\]
Indeed, the radial flux is the integral of the right side over a ball;
it vanishes outside its support. Set the exterior constant equal to zero.

Test \eqref{eq:linear-error} with $\psi_R$ and its first derivatives,
using $\ell_{CR}$ in that equation.
Integration by parts, using the bounds
$|D^3\psi_R|+|\nabla(\rho_{2R}-\rho_R)|\le CR^{-d-1}$
and
$|D^2\psi_R|+|\rho_{2R}-\rho_R|\le CR^{-d}$
together with~\eqref{eq:local-current-mass}
and~\eqref{eq:periodic-L1-potential}, gives
\[
 |b(2R)-b(R)|\le C/R,\qquad |a(2R)-a(R)|\le C.
\]
The same proof applies to $tR$ in place of $2R$, uniformly for
$1\le t\le2$. The gradient increments are summable along dyadic
scales, while the mean increments accumulate at most logarithmically.
Consequently,
\begin{equation}\label{eq:affine-limit}
 b(R)\longrightarrow b,\qquad |b(R)-b|\le C/R,
 \qquad |a(R)|\le C\log(2+R).
\end{equation}

To obtain the energy bounds for $g_0=h-b\cdot x$, write $\ell_R=\beta_R\cdot x+c_R$. The case $j=0$ of~\eqref{eq-mixed-energy} and the Cauchy--Schwarz inequality give $|b(R)-\beta_{CR}|\le C$. Hence \eqref{eq:affine-limit} bounds
$|\beta_{CR}-b|$ uniformly. Adding this bounded constant gradient to
\eqref{eq-mixed-energy}, and using the mixed trace mass bound there,
proves \eqref{eq:fixed-gradient} and \eqref{eq:fixed-current}.

Poincar\'e's inequality, with the weighted mean $\rho_{CR}*g_0(0)$, gives
$\fint_{B_R}|g_0-\rho_{CR}*g_0(0)|^2\le CR^2$.
Since $\rho$ is even, this weighted mean equals $a(CR)$.
Equation \eqref{eq:affine-limit} now proves \eqref{eq:fixed-energy}.
\end{proof}

\begin{proposition}
The vector $b$ in Lemma~\ref{lem:affine} satisfies
\begin{equation}\label{eq:averaged-log}
 \fint_{B_R}|u-Q-b\cdot x|\,dV\le C\log(2+R)
 \quad(R\ge2).
\end{equation}
Consequently,
\begin{equation}\label{eq:averaged-log-mu}
 \int_{B_R}|u-Q-b\cdot x|\,d\mu\le CR^d\log(2+R)
 \quad(R\ge2).
\end{equation}
\end{proposition}
\begin{proof}

We first represent $g_0-p$ as the sum of a renormalized singular integral of $S:=S(g_0)$ and an entire $A$-harmonic function. We show that the harmonic term is constant. The weak $(1,1)$ estimate and~\eqref{eq:fixed-energy} then yield~\eqref{eq:averaged-log}.

Since $S$ is positive, the total variation of each coefficient measure is bounded by a dimensional multiple of $S\wedge\omega_0$. Let $G_A$ be the fundamental solution of $L_A$, with $L_AG_A=\delta_0$. In real coordinates, write $dd^cS=\sum_\alpha D_\alpha S_\alpha$, where $S_\alpha$ are the coefficient measures of $S$ and $D_\alpha$ are constant-coefficient second-order operators. For a smooth compactly supported coefficient vector $f=(f_\alpha)$, the operator $T=L_A^{-1}dd^c$ has the representation
\[
 \begin{aligned}
 Tf(x)
 &=\sum_\alpha\bigl((D_\alpha G_A)*f_\alpha\bigr)(x)\\
 &=\sum_\alpha\left[
 c_\alpha f_\alpha(x)
 +\operatorname{p.v.}\int_{\R^d}K_\alpha(x-y)f_\alpha(y)\,dy
 \right].
 \end{aligned}
\]
Here $c_\alpha$ are constants and $K_\alpha(x)=D_\alpha G_A(x)$ for $x\ne0$. Writing $K=(K_\alpha)_\alpha$, we have
\begin{equation}\label{eq:CZ-kernel}
 |K(x)|\le C|x|^{-d},\qquad |\nabla K(x)|\le C|x|^{-d-1}
 \quad(x\ne0).
\end{equation}
The operators in the principal-value terms are $L^2$-bounded, since their Fourier multipliers are bounded. By the classical Calder\'on--Zygmund weak $(1,1)$ estimate (see~\cite[Chapter~II, \S\S2--3]{Stein1970} and~\cite[Chapter~I, \S5]{Stein1993}) and the $L^1$ boundedness of the terms $c_\alpha f_\alpha$,
\[
 \bigl|\{x\in\R^d:|Tf(x)|>t\}\bigr|
 \le\frac{C}{t}\sum_\alpha\norm{f_\alpha}_{L^1(\R^d)}
 \qquad(t>0)
\]
for every $L^1$ coefficient vector $f$.

We define the renormalized distribution $F=TS$ by replacing $K(x-y)$ with $K(x-y)-K(-y)$ for $|y|>1$. By~\eqref{eq:fixed-current} and~\eqref{eq:CZ-kernel}, the resulting tail integral converges at infinity, uniformly for $x$ in compact sets. For each $R\ge2$, this definition gives on $B_R$
\begin{equation}\label{eq:CZ-split}
 F=T(\mathbf1_{B_{2R}}S)+c_R+F_R.
\end{equation}
Summing $O(1)$ over dyadic annuli up to $2R$ gives
\begin{equation}\label{eq:CZ-constant-bound}
 |c_R|\le C\log(2+R).
\end{equation}
For the tail, summing $CR/r$ over annuli of radius $r\ge2R$,
using the gradient bound in \eqref{eq:CZ-kernel}, gives
\[
 \norm{F_R}_{L^\infty(B_R)}\le C.
\]

Equation \eqref{eq:linear-error}, with $\ell=b\cdot x$, implies that
\begin{equation}\label{eq:harmonic-remainder}
 H=g_0-p-F
\end{equation}
is an entire $A$-harmonic distribution. Elliptic regularity makes $H$
smooth; hence $F$ is locally integrable, since $g_0$ is continuous and
$p$ is locally integrable. By mollifying the local term in~\eqref{eq:CZ-split}, applying the weak $(1,1)$ estimate, and passing to the limit almost everywhere, we obtain
\begin{equation}\label{eq:weak-F}
 |\{x\in B_R:|F(x)-c_R|>t\}|
 \le CR^d\min\{1,t^{-1}\}\quad(t>0).
\end{equation}
The local integrability established above justifies this limiting step, including the case where $T$ contains an identity component.

Fix $0<q<1$. Equations \eqref{eq:fixed-energy},
\eqref{eq:CZ-constant-bound}, \eqref{eq:weak-F}, and \eqref{eq:periodic-L1-potential} imply
$(\fint_{B_R}|H|^q)^{1/q}\le CR$.
Interior estimates for harmonic functions then give
$|D^2H(x)|\le C_xR^{-1}$ on each fixed ball as $R\to\infty$.
Thus $H$ is affine.

To eliminate the linear part of $H$, we use the choice of $b$:
$\nabla(\rho_R*g_0)(0)=b(R)-b\to0$ by~\eqref{eq:affine-limit}.
Let $\kappa_R=T^*(\nabla\rho_R)$, where $T^*$ is the formal adjoint of $T$, applied to each component of $\nabla\rho_R$. Homogeneity gives
\[
 \kappa_R(y)=R^{-d-1}\kappa_1(y/R).
\]
The bounded Fourier multipliers of $T^*$ and the smoothness of $\rho$ imply that $\kappa_1$ is locally bounded. At infinity, the cancellation $\int\nabla\rho=0$ and the kernel bounds in~\eqref{eq:CZ-kernel} give the decay estimate $|\kappa_1(y)|\le C(1+|y|)^{-d-1}$. Integrating against the coefficient measures of $S$ and summing over dyadic annuli using~\eqref{eq:fixed-current} yields
\begin{equation}\label{eq:CZ-gradient}
 |\nabla(\rho_R*F)(0)|
 \le C\int (R+|y|)^{-d-1}\,d|S|(y)\le C/R.
\end{equation}
This is obtained first for truncated measures and then by the
renormalized limit; the subtracted constants disappear under the
gradient. Also \eqref{eq:periodic-L1-potential} gives
$|\nabla(\rho_R*p)(0)|\le C/R$.
Taking the smoothed gradient of \eqref{eq:harmonic-remainder} and using
\eqref{eq:affine-limit} and \eqref{eq:CZ-gradient} proves $\nabla H=0$.

By \eqref{eq:periodic-L1-potential}, Chebyshev's inequality gives
$|\{|p|>t\}\cap B_R|\le CR^d/t$. Adding this estimate to
\eqref{eq:weak-F}, and absorbing the constant $H$, proves the same
weak $L^1$ bound for $g_0-c_R$.
To make the tail integrable, we also use the coarse $L^2$ bound
from~\eqref{eq:fixed-energy}:
$\fint_{B_R}|g_0-c_R|^2\le CR^2$.
Integrating the smaller of these two distribution-function bounds gives
\[
 \fint_{B_R}|g_0-c_R|
 \le C\int_0^\infty\min\{1,t^{-1},R^2t^{-2}\}\,dt
 \le C\log(2+R).
\]
Combine this with \eqref{eq:CZ-constant-bound} to obtain
\eqref{eq:averaged-log}.

Fix $0\le\chi\le1$ in $C_c^\infty(B_2)$ with $\chi=1$ on $B_1$, and put $\chi_R(x)=\chi(x/R)$. Using the periodic corrector $\phi_A$ from Lemma~\ref{lem:measure-corrector}, set on $B_4$
\[
 q_R(z)=Q(z)+R^{-1}b\cdot z,\qquad
 a_R(z)=Q(z)+R^{-2}\phi_A(Rz).
\]
The growth bounds make $u_R,q_R,a_R,Q$ uniformly bounded psh functions, with $\MA(a_R)=\mu_R$, $\MA(Q)=dV$, and $\norm{a_R-Q}_{L^\infty(B_4)}\le CR^{-2}$. Apply Lemma~\ref{lem:local-transfer} to $(u_R,q_R,a_R,Q)$ and $(q_R,u_R,a_R,Q)$. Since $u_R-q_R=R^{-2}g_0(R\,\cdot)$, adding the estimates for the two positive parts and scaling back gives
\[
 \begin{aligned}
 \left|\int\chi_R|g_0|\,d\mu-\int\chi_R|g_0|\,dV\right|
 &=R^{d+2}\left|\int\chi\,|u_R-q_R|\,(d\mu_R-dV)\right|\\
 &\le CR^d.
 \end{aligned}
\]
Since $\chi_R=1$ on $B_R$, positivity of $\mu$ gives
\[
 \int_{B_R}|g_0|\,d\mu
 \le\int_{B_{2R}}|g_0|\,dV+CR^d.
\]
Equation~\eqref{eq:averaged-log} now proves~\eqref{eq:averaged-log-mu}.
\end{proof}

The next estimate supplies the closeness to $Q$ needed for the local
comparison argument in Section~\ref{sec:lowpoint}.

\begin{lemma}
For some $\gamma>0$, the vector $b$ of Lemma~\ref{lem:affine} satisfies
\begin{equation}\label{eq:flatness}
 \sup_{B_R}|u-Q-b\cdot x|\le CR^{2-\gamma}\qquad(R\ge2).
\end{equation}
\end{lemma}
\begin{proof}
Use the effective solution and quadratic polynomial from
Proposition~\ref{prop-quadratic-approx}, with a fixed larger radius.
The uniform homogenization estimate, followed by the interior estimate
for the effective solution minus that polynomial, gives
\[
 \sup_{B_R}|u-P_{2R}|\le CR^{2-\gamma}+CR.
\]
Here we take $\gamma=1/(2n+1)$, as given by Lemma~\ref{lem:uniform-homogenization}.
Equation~\eqref{eq-quadratic-coefficient-rate} gives
$\norm{\HR(P_{2R}-Q)}\le C/R$, while the constant-term estimate
in the proof of Proposition~\ref{prop-quadratic-approx} gives
$|P_{2R}(0)|\le CR$.
The proof of Lemma~\ref{lem:affine} bounds the difference between the
linear coefficient of $P_{2R}$ and $b$ by $C$.
Consequently $\sup_{B_R}|P_{2R}-Q-b\cdot x|\le CR$.
Adding this to the first displayed inequality proves \eqref{eq:flatness}.
\end{proof}

\section{Rigidity and completion of the main theorem}
\label{sec:lowpoint}

Inspired by the strategy of Jin--Li--Tran--Tu~\cite{JLT} for the real Monge--Amp\`ere equation, we combine the preceding estimates with an oscillation argument to complete the proof of Theorem~\ref{thm:main}.

Let $Q$ be the quadratic polynomial obtained in Proposition~\ref{prop-quadratic-approx}, and let $b\in\R^d$ be the vector from Lemma~\ref{lem:affine}.  
Let $\phi=\phi_A$, with $A=\Hess Q$, be the periodic corrector
satisfying~\eqref{eq:periodic-corrector-general}.
Define
\begin{equation}\label{eq:closure-g}
 \psi(x)=Q(x)+b\cdot x+\phi(x),\qquad
 g=u-\psi.
\end{equation}
Both $u$ and $\psi$ are continuous plurisubharmonic functions with
Monge--Amp\`ere measure $\mu$.
Since $\phi$ is bounded,  \eqref{eq:averaged-log-mu}  yields
\begin{equation}\label{eq:closure-weighted-log}
 \int_{B_R}|g|\,d\mu\le CR^d\log(2+R)\qquad(R\ge2).
\end{equation}

For $r\ge1$ and $y\in\R^d$, write
\[
 \mu_{r,y}(E)=\mu_r(E+y/r)
 \quad\text{for Borel sets }E.
\]

First, adapting the compactness and envelope arguments of
Jin--Li--Tran--Tu~\cite[Lemmas~2.7 and~2.10]{JLT} to the complex
setting, we establish the following large-scale estimates.

\begin{lemma}
\label{lem:closure-deficit}
Let $Q$ be a positive definite homogeneous real quadratic polynomial with $\mathcal M(Q)=dV$, and let $\mu$ be a nonnegative $\Gamma$-periodic measure satisfying~\eqref{eq:holder-measure}, with $\bar\mu=1$. Define
\[
 D=\{x:Q(x)<1\},\qquad K=\tfrac12\overline D.
\]
For every $\eta\in(0,1)$, there exist $\varepsilon\in(0,\tfrac12)$ and $r_0\ge1$, depending only on $\eta,Q,\Gamma,\mu$, such that the following holds. If $r\ge r_0$, $y\in\R^d$, and $q\in\operatorname{PSH}(D)\cap L^\infty(D)$ satisfy
\begin{equation}\label{eq:closure-q-hyp}
 \mathcal M(q)\le\mu_{r,y},\qquad Q-\varepsilon\le q\le1+\varepsilon,
 \qquad \inf_K(q-Q)\le\varepsilon,
\end{equation}
then
\[
 \mu_{r,y}(D)-\mathcal M(q)(D)<\eta\mu_{r,y}(D).
\]
\end{lemma}

\begin{proof}
Suppose the assertion fails. Then for some fixed $\eta\in(0,1)$ there are $\varepsilon_j\to0$, $r_j\to\infty$, phases $y_j$, and functions $q_j$ satisfying~\eqref{eq:closure-q-hyp}, with
\[
 \mathcal M(q_j)(D)\le(1-\eta)\mu_{r_j,y_j}(D).
\]
Let $\widehat q_j$ solve
\[
 \mathcal M(\widehat q_j)=\mathcal M(q_j)\quad\hbox{in }D,
 \qquad \widehat q_j=1\quad\hbox{on }\partial D.
\]

Let $\phi=\phi_{\Hess Q}$ be the periodic corrector from Lemma~\ref{lem:measure-corrector}. The functions
\[
 a_{r,y}(x)=Q(x)+r^{-2}\phi(y+rx),\qquad
 \mathcal M(a_{r,y})=\mu_{r,y},
\]
are uniformly bounded in $C^{0,\alpha}$ on a fixed neighborhood of $\overline D$ for some $\alpha\in(0,1]$, since the periodic part has H\"older seminorm at most $r^{\alpha-2}[\phi]_{C^{0,\alpha}}$. Let $H_{r,y}$ be the maximal psh function with boundary value $-a_{r,y}$. Since $D$ is smooth and strongly pseudoconvex, the homogeneous Dirichlet estimates~\cite{BT76} give uniform H\"older bounds for $H_{r,y}$. Thus $\zeta_{r,y}=a_{r,y}+H_{r,y}$ satisfies
\[
 \mathcal M(\zeta_{r,y})\ge\mu_{r,y},\qquad
 \zeta_{r,y}|_{\partial D}=0,\qquad
 \norm{\zeta_{r,y}}_{C^{0,\alpha_*}(\overline D)}\le C
\]
for some fixed $\alpha_*>0$. Comparison~\cite[Theorem~4.1]{BT} gives $1+\zeta_{r_j,y_j}\le\widehat q_j\le1$. Nguyen's Dirichlet regularity theorem~\cite[Theorem B and its proof in Section~4]{NguyenII} therefore yields
\[
 \norm{\widehat q_j}_{C^{0,\beta}(\overline D)}\le C
\]
for some $\beta>0$ and $C$ independent of $j$.

Since $Q=1$ on $\partial D$, the bounds in~\eqref{eq:closure-q-hyp} place the boundary liminf and limsup of $q_j$ between $1-\varepsilon_j$ and $1+\varepsilon_j$. Comparison with $\widehat q_j\pm\varepsilon_j$ gives
\[
 \norm{q_j-\widehat q_j}_{L^\infty(D)}\le\varepsilon_j.
\]
Passing to a subsequence, we obtain uniform convergence of $q_j$ on $D$ to a function $q\in\operatorname{PSH}(D)\cap C(\overline D)$ satisfying $q\ge Q$ and $q=1$ on $\partial D$.

Since the period converges to zero and $\bar{\mu}=1$, we have $\mu_{r_j,y_j}\rightharpoonup dV$ and
\[
 \mu_{r_j,y_j}(D)=|D|+O(r_j^{-1})\to|D|,
\]
uniformly in the phase. Bedford--Taylor continuity~\cite{BT} and a compact exhaustion of $D$ then give
\[
 \mathcal M(q)\le dV,\qquad
 \mathcal M(q)(D)\le\liminf_{j\to\infty}\mathcal M(q_j)(D)
                  \le(1-\eta)|D|.
\]
The last condition in~\eqref{eq:closure-q-hyp}, together with uniform convergence on the compact set $K\Subset D$, yields a point $x_0\in K$ such that $q(x_0)=Q(x_0)$. By~\cite[Proposition~3.11(1)]{Zeriahi2013}, $q$ is a viscosity supersolution of $\det\Hess q=1$. On the cone of real Hessians with positive definite complex part, the complex determinant is locally strictly elliptic and locally Lipschitz, and extends continuously with value zero to the boundary of this cone. Since $Q$ is smooth and strictly plurisubharmonic, the strong comparison principle of Li--Nguyen--Wang~\cite[Theorem~2.3]{LiNguyenWang2020} gives $q=Q$; see also Caffarelli--Li--Nirenberg~\cite{CLN2013}. This contradicts the preceding mass bound.
\end{proof}

Let $\varepsilon,r_0$ be given by Lemma~\ref{lem:closure-deficit} for $\eta=1/4$, and set $\delta_0=\varepsilon/2$.

\begin{lemma}\label{lem:closure-lowpoint}
There is $M_0\ge1$, independent of $r,y,s,c$, such that the following holds. Suppose $r\ge r_0$, $y\in\R^d$, and $U\in\operatorname{PSH}(D)\cap C(\overline D)$ satisfies
\begin{equation}\label{eq:closure-U-hyp}
 \mathcal M(U)=\mu_{r,y},\qquad
 \norm{U-Q}_{L^\infty(D)}\le\delta_0.
\end{equation}
Let $V\in\operatorname{PSH}(D)\cap C(\overline D)$ satisfy $\mathcal M(V)=\mu_{r,y}$.
For any $s>0$ and $c\in\R$, set
\begin{equation}\label{eq:closure-v}
 v=\frac{U-V+c}{s}.
\end{equation}
If $v\ge0$ on $D$ and $\inf_Kv\le1$, then
\begin{equation}\label{eq:closure-lowpoint}
 \mu_{r,y}\{x\in D:v(x)\le M_0\}
                 \ge\tfrac34\mu_{r,y}(D).
\end{equation}
\end{lemma}

\begin{proof}
Choose $L\ge1$ so that $L^{-1}\le\varepsilon/2$, and let
\[
 h=\min\{U+v/L,\,1+\delta_0\},\qquad
 q=P_D(h):=\bigl(\sup\{p\in\operatorname{PSH}(D):p\le h\}\bigr)^*,
\]
where $*$ denotes upper semicontinuous regularization. Since $h$ is continuous and $U\le h$, the envelope is bounded and satisfies $U\le q\le h$. On a ball compactly contained in $\{q<h\}$, a sufficiently small convex combination with a maximal psh modification remains below $h$; the envelope property therefore makes $q$ maximal there. Thus $\mathcal M(q)=0$ on $\{q<h\}$. The Bedford--Taylor contact inequality (see~\cite[Lemma~1.3]{GLZ}), applied to $q\le1+\delta_0$, also gives $\mathcal M(q)=0$ on $\{q=1+\delta_0\}$. Hence $\mathcal M(q)$ is carried by
\[
 E=\{q=U+v/L\},\qquad t=Ls>0.
\]

By~\eqref{eq:closure-v},
\[
 V+tq\le(1+t)U+c\quad\hbox{on }D,
 \qquad V+tq=(1+t)U+c\quad\hbox{on }E.
\]
Both sides are psh. The contact inequality and positivity of the mixed Monge--Amp\`ere products give
\begin{equation}\label{eq:closure-contact-upper}
 t^n{\bf1}_E\mathcal M(q)
 \le{\bf1}_E\mathcal M(V+tq)
 \le(1+t)^n{\bf1}_E\mu_{r,y}.
\end{equation}
It follows that $\mathcal M(q)=a\,d\mu_{r,y}$, with $a=0$ $\mu_{r,y}$-almost everywhere off $E$. Since $\mu_{r,y}=\mathcal M(U)$ does not charge pluripolar sets, Dinew's mixed Monge--Amp\`ere inequality~\cite{Dinew} and~\eqref{eq:closure-contact-upper} yield
\[
 \bigl(1+t a^{1/n}\bigr)^n\le(1+t)^n
 \qquad\mu_{r,y}\text{-almost everywhere on }E.
\]
Thus $a\le1$, and
\begin{equation}\label{eq:closure-envelope-measure}
 \mathcal M(q)\le\mu_{r,y},\qquad
 \mathcal M(q)(D\setminus E)=0.
\end{equation}

The bounds on $U$ and $\inf_Kv\le1$ give
\[
 Q-\varepsilon\le q\le1+\varepsilon,\qquad
 \inf_K(q-Q)\le\delta_0+L^{-1}\le\varepsilon.
\]
Lemma~\ref{lem:closure-deficit} and~\eqref{eq:closure-envelope-measure} therefore imply
\[
 \mu_{r,y}(E)\ge\mathcal M(q)(D)>\tfrac34\mu_{r,y}(D).
\]
On $E$, we have
\[
 v=L(q-U)\le L(1+2\delta_0)\le2L,
\]
so $M_0=2L$ gives~\eqref{eq:closure-lowpoint}.
\end{proof}

For $v=u$ or $v=\psi$, subtract the affine function
\[
 T_y(x)=Q(y)+DQ(y)\cdot(x-y)+b\cdot x
\]
and set $v_{r,y}(z)=r^{-2}[v(y+rz)-T_y(y+rz)]$.
The homogeneity of $Q$ gives
\begin{equation}\label{eq:closure-scaled-flatness}
 v_{r,y}(z)-Q(z)
   =r^{-2}\bigl[v(y+rz)-Q(y+rz)-b\cdot(y+rz)\bigr],
 \qquad \mathcal M(v_{r,y})=\mu_{r,y}.
\end{equation}
Choose $0<\gamma<2$ as in \eqref{eq:flatness}.
For $y+rD\subset B_{2R}$, that estimate and the boundedness of
$\phi$ show that the first expression in
\eqref{eq:closure-scaled-flatness} has absolute value at most
$CR^{2-\gamma}/r^2$, for both choices of $v$.
For a fixed sufficiently large $C_*$, the conditions
\begin{equation}\label{eq:closure-minscale}
 r_*(R):=C_*R^{1-\gamma/2},\qquad
 r\ge r_*(R),\quad y+rD\subset B_{2R},
\end{equation}
ensure that $v_{r,y}$ satisfies~\eqref{eq:closure-U-hyp} for both choices of $v$, once $R$ is large. Increase $C_*$ if necessary to include the
fixed lower bound $r_0$.
After applying the same rescaling to the other solution, we see that any nonnegative function of the form
\begin{equation}\label{eq:closure-physical-difference}
 \frac{u-\psi+c}{s}\quad\hbox{or}\quad\frac{\psi-u+c}{s},
 \qquad s>0,
\end{equation}
is precisely of the form \eqref{eq:closure-v}, with a rescaled
positive denominator. Thus, whenever \eqref{eq:closure-minscale}
holds, Lemma~\ref{lem:closure-lowpoint} applies on $y+rD$ to either
difference and every $s>0$, provided its infimum on $y+rK$
is at most one. Moreover, $r_*(R)/R=C_*R^{-\gamma/2}\to0$.

In the real case, the dichotomous Harnack argument of Jin--Li--Tran--Tu~\cite[Theorem~2.14]{JLT} yields either a local upper bound or exponential growth along a chain of points. Using the averaged logarithmic bound and Lemma~\ref{lem:closure-lowpoint}, we adapt this argument and rule out the exponential growth alternative. This gives the following pointwise logarithmic estimate.

\begin{proposition}
For all $R\ge2$,
\begin{equation}\label{eq:closure-pointwise-log}
 \sup_{B_R}|g|\le C\log(2+R).
\end{equation}
\end{proposition}

\begin{proof}
For a large fixed $R$, take either sign in
$G=\pm g/\log(2+R)$. The weighted logarithmic estimate
\eqref{eq:closure-weighted-log} gives a constant $C_1$ such that
\begin{equation}\label{eq:closure-tail}
 \mu\bigl(\{|G|>s\}\cap B_{2R}\bigr)
                  \le C_1R^d/s\qquad(s>0).
\end{equation}
The original quadratic bounds on $u$ and the definition of $\psi$ also
give the crude bound
\begin{equation}\label{eq:closure-crude}
 \sup_{B_{2R}}|G|\le C_2R^2.
\end{equation}

Put
\[
 \nu_+=1+\frac1{2M_0},\qquad
 \nu=1+\frac1{4M_0},\qquad 1<\nu<\nu_+.
\]
For all sufficiently large $r$, periodic cell counting gives
\begin{equation}\label{eq:closure-mass-comparable}
 c_0r^d\le\mu(y+rD)\le C_0r^d,
\end{equation}
with $c_0,C_0>0$ independent of $y$. Choose a fixed $L_1\ge1$ so large that $\tfrac34c_0L_1^d>2C_1$.
At a point $y$ with $G(y)=s>0$, consider
\begin{equation}\label{eq:closure-chain-radius}
 r=L_1\max\{Rs^{-1/d},\,r_*(R)\}.
\end{equation}
Suppose $y+rD\subset B_{2R}$. We claim that
\begin{equation}\label{eq:closure-chain-growth}
 \sup_{y+rD}G>\nu_+s.
\end{equation}
Otherwise the nonnegative function
\[
 v=\frac{\nu_+s-G}{(\nu_+-1)s}
\]
has $v(y)=1$ and has the form
\eqref{eq:closure-physical-difference}.
By~\eqref{eq:closure-minscale} and Lemma~\ref{lem:closure-lowpoint}, the set in $y+rD$ where
\[
 G\ge\bigl[\nu_+-M_0(\nu_+-1)\bigr]s
      =\left(\tfrac12+\tfrac1{2M_0}\right)s>s/2
\]
has $\mu$-measure at least $\tfrac34\mu(y+rD)$. This measure is at least $\tfrac34c_0r^d$, by
\eqref{eq:closure-mass-comparable}, whereas
\eqref{eq:closure-tail} bounds it by $2C_1R^d/s$.
These inequalities contradict \eqref{eq:closure-chain-radius} and
the choice of $L_1$. This proves \eqref{eq:closure-chain-growth}.


Let $C_D=\sup_{x\in D}|x|$.
Choose a fixed $M\ge1$ so large that
\[
 \frac{C_DL_1M^{-1/d}}{1-\nu^{-1/d}}<\frac1{16}.
\]
Set
\begin{equation}\label{eq:closure-chain-choices}
 N=\left\lfloor\frac{R}{16C_DL_1r_*(R)}\right\rfloor.
\end{equation}
For sufficiently large $R$, $N\ge cR^{\gamma/2}$ and
$C_DL_1r_*(R)<R/16$.
If some $y_0\in B_{R/2}$ satisfies $s_0=G(y_0)>M$,
apply \eqref{eq:closure-chain-growth} repeatedly and choose
$y_{j+1}\in y_j+r_jD$ with
\begin{equation}\label{eq:closure-chain-recursion}
 s_{j+1}:=G(y_{j+1})\ge\nu s_j,
 \qquad
 r_j=L_1\max\{Rs_j^{-1/d},r_*(R)\}.
\end{equation}
For every $k\le N$, the total displacement is at most
\[
 \sum_{j=0}^{k-1}|y_{j+1}-y_j|
 \le C_DL_1\left(
       R M^{-1/d}\sum_{j\ge0}\nu^{-j/d}+Nr_*(R)\right)
       <R/8.
\]
Each ellipsoid $y_j+r_jD$ is contained in $B_{C_Dr_j}(y_j)$, where $C_Dr_j\le C_DL_1(RM^{-1/d}+r_*(R))<R/8$ for sufficiently large $R$.
Thus all centers lie in $B_{5R/8}$ and all ellipsoids lie in
$B_{3R/4}\subset B_{2R}$; the hypothesis used in
\eqref{eq:closure-chain-growth} is justified at every step.
Equations \eqref{eq:closure-chain-recursion} and
\eqref{eq:closure-chain-choices} now imply
\[
 G(y_N)\ge M\nu^N\ge M\exp(cR^{\gamma/2}),
\]
contradicting \eqref{eq:closure-crude} for large $R$.
Therefore $G\le M$ on $B_{R/2}$.
Applying the argument to both signs and replacing $R$ by $2R$
proves \eqref{eq:closure-pointwise-log}; bounded radii are absorbed
into the constant by continuity.
\end{proof}

\begin{proposition}
There are $\varepsilon\in(0,1)$ and $R_0$ such that
\begin{equation}\label{eq:closure-contraction}
 \osc_{\tfrac{R}{2}D}g\le(1-\varepsilon)\osc_{RD}g
 \qquad(R\ge R_0).
\end{equation}
\end{proposition}

\begin{proof}
On $RD$, put $m_R=\inf g$, $M_R=\sup g$, and
$O_R=M_R-m_R$. If $O_R=0$, there is nothing to prove.
Otherwise set $v=(g-m_R)/O_R\in[0,1]$.
At least one of the sets $\{v\ge1/2\}$ and
$\{1-v\ge1/2\}$ has $\mu$-measure at least $\mu(RD)/2$.

For sufficiently large $R$, the rescaled functions $u_{R,0}$ and $\psi_{R,0}$ both satisfy~\eqref{eq:closure-U-hyp}, by~\eqref{eq:flatness} and~\eqref{eq:closure-scaled-flatness}.
Suppose first that $\mu\{v\ge1/2\}\ge\mu(RD)/2$.
If $\inf_{\tfrac{R}{2}D}(4M_0v)\le1$, then after rescaling we apply Lemma~\ref{lem:closure-lowpoint} with $U=u_{R,0}$. The estimate~\eqref{eq:closure-lowpoint} gives
\[
 \mu\{v\le1/4\}\ge\tfrac34\mu(RD),
\]
a contradiction. Hence $v\ge1/(4M_0)$ on $\tfrac{R}{2}D$.
In the other case, take $U=\psi_{R,0}$ after rescaling and apply the same argument to $4M_0(1-v)$; it gives
$v\le1-1/(4M_0)$ on $\tfrac{R}{2}D$.
In either case, \eqref{eq:closure-contraction} follows with
$\varepsilon=1/(4M_0)$.
\end{proof}

\begin{proof}[\normalfont\bfseries Proof of Theorem~\ref{thm:main}]

Fix $r\ge R_0$. Iterating \eqref{eq:closure-contraction} from
the outer ellipsoid $2^krD$ down to $rD$, and then using
\eqref{eq:closure-pointwise-log}, gives
\[
 \osc_{rD}g
 \le(1-\varepsilon)^k\osc_{2^krD}g
 \le C(1-\varepsilon)^k
           \log(2+C_D2^kr)\longrightarrow0.
\]
Thus $g$ is constant on every sufficiently large ellipsoid, hence
on $\C^n$. By \eqref{eq:closure-g}, there is $c\in\R$ such that
\begin{equation}\label{eq:closure-final-decomposition}
 u(x)=Q(x)+b\cdot x+c+\phi(x).
\end{equation}
The corrector is periodic, H\"older continuous, and normalized by
\eqref{eq:periodic-corrector-general}; the quadratic bounds and determinant
identity for $Q$ were proved in
Proposition~\ref{prop-quadratic-approx}. This proves existence in the normalized case. Multiplying the decomposition by the original mean mass to the power $1/n$ restores~\eqref{eq:decomposition}--\eqref{eq:matrix-main} for general $\mu$.

The quadratic part $Q$ is uniquely determined by the blow-down limit.
Lemma~\ref{lem:measure-corrector} then determines the normalized periodic
part $\phi$, and hence also the affine part $\ell=u-Q-\phi$.
\end{proof}

\appendix
\section{Proof of the smooth finite-ball estimate}\label{app:estimate}

We prove the smooth case of Proposition~\ref{prop:finite}. Throughout this appendix, write $m=m_0$ and $M=M_0$, and set
\[
A=\Hess w,\qquad
\mathcal L_w b=\tr(A^{-1}\Hess b),\qquad S=\tr A^{-1}.
\]
Thus $\mathcal L_w w=n$ and $\mathcal L_w |z|^2=S$.
Fix $0<t\le1$, and define
\[
Q_t=\G_t w-w,\qquad Z_t=\Sph_t w-w,
\]
where
\[
 \Sph_t w(z)=\frac1{|\partial B_t|}
 \int_{\partial B_t}w(z+h)\,d\sigma(h)
\]
denotes the spherical mean centered at $z$, and $d\sigma$ is surface measure.

The divergence theorem, applied to the spherical mean, gives
\[
\frac{d}{ds}\Sph_s w
=\frac{s}{2n}\B_s(\Delta_{\R}w)
=\frac{2s}{n}\B_s(\tr\Hess w).
\]
Integrating and using \eqref{eq:G} yields
\begin{equation}\label{eq:appendix-average}
Z_t=\frac{t^2}{n}\tr\Hess(\G_t w).
\end{equation}
Monotonicity of spherical means and plurisubharmonicity also give
\begin{equation}\label{eq:appendix-average-positivity}
0\le Q_t\le Z_t,\qquad \Hess(\G_t w)\ge0.
\end{equation}

At a fixed point $z$, the matrices $A(z+h)$ and $A(z)$ have the same
determinant. The arithmetic--geometric mean inequality gives
$\tr[A(z)^{-1}A(z+h)]\ge n$.
Average this inequality with the kernel of $\G_t$ to obtain
\begin{equation}\label{eq:appendix-source}
\mathcal L_w Q_t\ge0.
\end{equation}
The coefficients of $\mathcal L_w $ stay fixed at $z$ during this average.

Introduce the positive comparison function
\begin{equation}\label{eq:appendix-V}
V=w-\frac m2|z|^2,\qquad
\frac m2\rho\le V\le M\rho,\qquad
\mathcal L_w V=n-\frac m2S.
\end{equation}
Set $T=4n/m$ and $\lambda_*=T^{1-n}$, as in
\eqref{eq:finite-constants}. Where $S\le T$, every eigenvalue of
$A^{-1}$ is at most $T$, whereas their product is one.
Therefore $A^{-1}\ge\lambda_*\Id$ there.

Using this lower matrix bound and
\eqref{eq:appendix-average}--\eqref{eq:appendix-average-positivity} gives
\[
Q_t\le\frac{t^2}{n}\tr\Hess(\G_t w)
\le\frac{t^2}{n\lambda_*}\mathcal L_w (\G_t w).
\]
Since $\mathcal L_w w=n$, we have proved
\begin{equation}\label{eq:appendix-good}
Q_t\le\frac{t^2}{\lambda_*}+\frac{t^2}{n\lambda_*}\mathcal L_w Q_t
\quad\hbox{where }S\le T.
\end{equation}

Take $t\le t_0=\min\{1,\sqrt{m\lambda_*/2}\}$, so that
$t^2/\lambda_*\le m/2\le V$. We choose the following quotient so that any positive interior maximum must lie in the region $S\le T$, where~\eqref{eq:appendix-good} applies:

\begin{equation}\label{eq:appendix-J}
H=Q_t-\frac{t^2}{\lambda_*},\quad
\beta=\frac{m}{16M+m},\quad a=1-\beta,\quad
J=\frac{H}{V^a\rho^{2\beta}}.
\end{equation}
The parameters satisfy $a+\beta=1$ and $4\beta=am/(4M)$.

Suppose $J$ has a positive interior maximum on $\overline B_L$.
At that point $H>0$ and
\[
\partial\log H=a\,\partial\log V+2\beta\,\partial\log\rho.
\]
Weighted Cauchy--Schwarz, with weights $a,\beta$, gives
\[
|\partial\log H|_{A^{-1}}^2
\le a|\partial\log V|_{A^{-1}}^2
+4\beta|\partial\log\rho|_{A^{-1}}^2.
\]
Apply the chain rule to $\mathcal L_w \log J$, and use this inequality.
Since $|\partial\log\rho|_{A^{-1}}^2\le S/\rho$, we obtain
\[
0\ge\mathcal L_w \log J
\ge\frac{\mathcal L_w H}{H}
-a\frac{\mathcal L_w V}{V}-\frac{4\beta S}{\rho}.
\]
By $V\le M\rho$ and $4\beta=am/(4M)$,
$4\beta S/\rho\le amS/(4V)$. Substituting
$\mathcal L_w V=n-mS/2$ from \eqref{eq:appendix-V} yields
\begin{equation}\label{eq:appendix-max}
\frac{\mathcal L_w H}{H}\le\frac aV\left(n-\frac m4S\right).
\end{equation}

By \eqref{eq:appendix-source}, $\mathcal L_w H=\mathcal L_w Q_t\ge0$.
Thus \eqref{eq:appendix-max} forces $S\le4n/m=T$.
We may therefore use \eqref{eq:appendix-good}, which gives
\[
H\le\frac{t^2}{n\lambda_*}\mathcal L_w H.
\]
As $H>0$ and $t^2/\lambda_*\le V$, the last inequality implies
$\mathcal L_w H/H\ge n/V$. This contradicts
\eqref{eq:appendix-max}, whose right side is at most
$an/V<n/V$. Hence $J$ has no positive interior maximum.

For $z\in B_L$, all translates in $\G_t$ lie in $B_{2L}$.
The upper quadratic bound and positivity of $w$ give
\[
H\le Q_t\le\G_t w\le M\G_t\rho
=M\left(\rho+\frac{nt^2}{2(n+1)}\right)\le2M\rho.
\]
On $\partial B_L$, \eqref{eq:appendix-V} and
\eqref{eq:appendix-J} consequently give
\[
J\le2M\left(\frac2m\right)^a(1+L^2)^{-\beta}.
\]
The maximum principle just proved extends this bound to $B_L$.
Multiplying by $V^a\rho^{2\beta}$ and using $V\le M\rho$ gives
\[
H(z)\le2M\left(\frac{2M}{m}\right)^a
\rho(z)\left(\frac{\rho(z)}{1+L^2}\right)^\beta.
\]
For $|z|<\sqrt L$, the last ratio is at most $2/L$.
Since $a+\beta=1$ and $M/m\ge1$, it follows that
\[
H(z)\le\frac{4M^2}{m}L^{-\beta}\rho(z).
\]
Substituting $Q_t=H+t^2/\lambda_*$ from~\eqref{eq:appendix-J} and using $Q_t\ge0$ from~\eqref{eq:appendix-average-positivity} gives~\eqref{eq:finite}.\hfill$\square$

\section*{Declaration of generative AI use}

During the preparation of this work, the authors used OpenAI's Codex to assist with literature searches, the discussion and checking of mathematical arguments, and manuscript drafting and revision, including \LaTeX{} editing. The authors take full responsibility for the mathematical results and the final content of the manuscript.

\begingroup
\color{black}
\hypersetup{linkcolor=black,citecolor=black,urlcolor=black}

\endgroup

\begingroup
\small
\hypersetup{urlcolor=black}
\bigskip

\noindent J. Hu

\nopagebreak[4]
\noindent School of Mathematics, Sichuan University\\
Chengdu, China\\[1mm]
Emails: \href{mailto:jingchenhoo@gmail.com}{\textsf{jingchenhoo@gmail.com}}; \href{mailto:jingchenhu@scu.edu.cn}{\textsf{jingchenhu@scu.edu.cn}}

\bigskip

\noindent X. Tu

\nopagebreak[4]
\noindent Department of Mathematics, The Hong Kong University of Science and Technology\\
Clear Water Bay, Kowloon, Hong Kong\\[1mm]
Emails: \href{mailto:tuxushan28561m@gmail.com}{\textsf{tuxushan28561m@gmail.com}}; \href{mailto:maxstu@ust.hk}{\textsf{maxstu@ust.hk}}
\par\endgroup

\begin{thebibliography}{45}\small

\bibitem{AL1989}
M.~Avellaneda and F.-H.~Lin,
Un th\'eor\`eme de Liouville pour des \'equations elliptiques \`a coefficients p\'eriodiques,
\emph{C. R. Acad. Sci. Paris S\'er. I Math.} \textbf{309} (1989), 245--250.

\bibitem{BCGJ2003}
J.~Bao, J.~Chen, B.~Guan, and M.~Ji,
Liouville property and regularity of a Hessian quotient equation,
\emph{Amer. J. Math.} \textbf{125} (2003), 301--316.
\href{https://doi.org/10.1353/ajm.2003.0007}{doi:10.1353/ajm.2003.0007}.

\bibitem{BT76}
E.~Bedford and B.~A.~Taylor, The Dirichlet problem for a complex
Monge--Amp\`ere equation, \emph{Invent. Math.} \textbf{37} (1976), 1--44.
\href{https://doi.org/10.1007/BF01418826}{doi:10.1007/BF01418826}.

\bibitem{BT}
E.~Bedford and B.~A.~Taylor,
A new capacity for plurisubharmonic functions,
\emph{Acta Math.} \textbf{149} (1982), 1--40.
\href{https://doi.org/10.1007/BF02392348}{doi:10.1007/BF02392348}.

\bibitem{BD}
Z.~B{\l}ocki and S.~Dinew,
A local regularity of the complex Monge--Amp\`ere equation,
\emph{Math. Ann.} \textbf{351} (2011), 411--416.
\href{https://doi.org/10.1007/s00208-010-0609-0}{doi:10.1007/s00208-010-0609-0}.

\bibitem{CKNS}
L.~A.~Caffarelli, J.~J.~Kohn, L.~Nirenberg, and J.~Spruck,
The Dirichlet problem for nonlinear second-order elliptic equations.
II. Complex Monge--Amp\`ere, and uniformly elliptic, equations,
\emph{Comm. Pure Appl. Math.} \textbf{38} (1985), 209--252.
\href{https://doi.org/10.1002/cpa.3160380206}{doi:10.1002/cpa.3160380206}.

\bibitem{CL2003}
L.~A.~Caffarelli and Y.~Y.~Li,
An extension to a theorem of J\"orgens, Calabi, and Pogorelov,
\emph{Comm. Pure Appl. Math.} \textbf{56} (2003), 549--583.
\href{https://doi.org/10.1002/cpa.10067}{doi:10.1002/cpa.10067}.

\bibitem{CL2004}
L.~A.~Caffarelli and Y.~Y.~Li,
A Liouville theorem for solutions of the Monge--Amp\`ere equation with periodic data,
\emph{Ann. Inst. H. Poincar\'e C Anal. Non Lin\'eaire} \textbf{21} (2004), 97--120.
\href{https://doi.org/10.1016/j.anihpc.2003.01.005}{doi:10.1016/j.anihpc.2003.01.005}.

\bibitem{CLN2013}
L.~A.~Caffarelli, Y.~Y.~Li, and L.~Nirenberg,
Some remarks on singular solutions of nonlinear elliptic equations III: viscosity solutions including parabolic operators,
\emph{Comm. Pure Appl. Math.} \textbf{66} (2013), 109--143.
\href{https://doi.org/10.1002/cpa.21412}{doi:10.1002/cpa.21412}.

\bibitem{Calabi}
E.~Calabi,
Improper affine hyperspheres of convex type and a generalization of a theorem by K.~J\"orgens,
\emph{Michigan Math. J.} \textbf{5} (1958), 105--126.
\href{https://doi.org/10.1307/mmj/1028998055}{doi:10.1307/mmj/1028998055}.

\bibitem{CHS2026}
H.~Chen, J.~Hu, and L.~Sheng,
Pogorelov interior estimates and a Liouville theorem for the complex Monge--Amp\`ere equation,
preprint, 2026.
\href{https://arxiv.org/abs/2609.18265}{arXiv:2609.18265}.

\bibitem{ChuLiu2026}
J.~Chu and Y.~Liu,
Quadratic growth rigidity for the complex Monge--Amp\`ere equation,
preprint, 2026.
\href{https://arxiv.org/abs/2610.26016}{arXiv:2610.26016}.

\bibitem{DDGHKZ}
J.-P.~Demailly, S.~Dinew, V.~Guedj, P.~H.~Hiep, S.~Ko{\l}odziej,
and A.~Zeriahi,
H\"older continuous solutions to Monge--Amp\`ere equations,
\emph{J. Eur. Math. Soc.} \textbf{16} (2014), 619--647.
\href{https://doi.org/10.4171/JEMS/442}{doi:10.4171/JEMS/442}.

\bibitem{Dinew}
S.~Dinew, An inequality for mixed Monge--Amp\`ere measures,
\emph{Math. Z.} \textbf{262} (2009), 1--15.
\href{https://doi.org/10.1007/s00209-008-0356-z}{doi:10.1007/s00209-008-0356-z}.

\bibitem{Evans}
L.~C.~Evans,
Classical solutions of fully nonlinear, convex, second-order elliptic
equations, \emph{Comm. Pure Appl. Math.} \textbf{35} (1982), 333--363.
\href{https://doi.org/10.1002/cpa.3160350303}{doi:10.1002/cpa.3160350303}.

\bibitem{GLZ}
V.~Guedj and C.~H.~Lu,
Quasi-plurisubharmonic envelopes 1: uniform estimates on K\"ahler manifolds,
\emph{J. Eur. Math. Soc.} \textbf{27} (2025), 1185--1208.
\href{https://doi.org/10.4171/JEMS/1460}{doi:10.4171/JEMS/1460}.

\bibitem{Hiep}
P.~H.~Hiep, H\"older continuity of solutions to the Monge--Amp\`ere
equations on compact K\"ahler manifolds,
\emph{Ann. Inst. Fourier} \textbf{60} (2010), 1857--1869.
\href{https://doi.org/10.5802/aif.2574}{doi:10.5802/aif.2574}.

\bibitem{HuSheng2026}
J.~Hu and L.~Sheng,
Interior $C^2$ estimates and a Liouville theorem for the complex Monge--Amp\`ere equation,
preprint, 2026.
\href{https://arxiv.org/abs/2610.28375}{arXiv:2610.28375}.

\bibitem{JLT}
T.~Jin, Y.~Y.~Li, H.~V.~Tran, and X.~Tu,
A Liouville theorem for convex functions with periodic Monge--Amp\`ere
measure, preprint, 2026.
\href{https://arxiv.org/abs/2511.15021v3}{arXiv:2511.15021v3}.

\bibitem{Jorgens}
K.~J\"orgens,
\"Uber die L\"osungen der Differentialgleichung $rt-s^2=1$,
\emph{Math. Ann.} \textbf{127} (1954), 130--134.
\href{https://doi.org/10.1007/BF01361114}{doi:10.1007/BF01361114}.

\bibitem{Kol1998}
S.~Ko{\l}odziej,
The complex Monge--Amp\`ere equation,
\emph{Acta Math.} \textbf{180} (1998), 69--117.
\href{https://doi.org/10.1007/BF02392879}{doi:10.1007/\allowbreak BF02392879}.

\bibitem{Kol2008}
S.~Ko{\l}odziej,
H\"older continuity of solutions to the complex Monge--Amp\`ere equation with the right-hand side in $L^p$: the case of compact K\"ahler manifolds,
\emph{Math. Ann.} \textbf{342} (2008), 379--386.
\href{https://doi.org/10.1007/s00208-008-0239-y}{doi:10.1007/s00208-008-0239-y}.

\bibitem{KN2019}
S.~Ko{\l}odziej and N.~C.~Nguyen,
H\"older continuous solutions of the Monge--Amp\`ere equation on compact
Hermitian manifolds, \emph{Ann. Inst. Fourier} \textbf{68} (2018), 2951--2964.
\href{https://doi.org/10.5802/aif.3232}{doi:10.5802/aif.3232}.

\bibitem{KNStability2019}
S.~Ko{\l}odziej and N.~C.~Nguyen,
Stability and regularity of solutions of the Monge--Amp\`ere equation on
Hermitian manifolds, \emph{Adv. Math.} \textbf{346} (2019), 264--304.
\href{https://doi.org/10.1016/j.aim.2019.02.004}{doi:10.1016/j.aim.2019.02.004}.

\bibitem{LiSheng2021}
A.-M.~Li and L.~Sheng,
A Liouville theorem on the PDE $\det(f_{i\bar j})=1$,
\emph{Math. Z.} \textbf{297} (2021), 1623--1632.
\href{https://doi.org/10.1007/s00209-020-02571-z}{doi:10.1007/s00209-020-02571-z}.

\bibitem{LiLiang2026}
D.~Li and L.~Liang,
Quadratic growth solutions of fully nonlinear elliptic equations with periodic data,
\emph{J. Lond. Math. Soc.} \textbf{113} (2026), e70446.
\href{https://doi.org/10.1112/jlms.70446}{doi:10.1112/jlms.70446}.

\bibitem{LRW2016}
M.~Li, C.~Ren, and Z.~Wang,
An interior estimate for convex solutions and a rigidity theorem,
\emph{J. Funct. Anal.} \textbf{270} (2016), 2691--2714.
\href{https://doi.org/10.1016/j.jfa.2016.01.008}{doi:10.1016/j.jfa.2016.01.008}.

\bibitem{LW2001}
P.~Li and J.~Wang,
Polynomial growth solutions of uniformly elliptic operators of non-divergence form,
\emph{Proc. Amer. Math. Soc.} \textbf{129} (2001), 3691--3699.
\href{https://doi.org/10.1090/S0002-9939-01-06167-6}{doi:10.1090/S0002-9939-01-06167-6}.

\bibitem{Li1990}
Y.~Y.~Li,
Some existence results for fully nonlinear elliptic equations of Monge--Amp\`ere type,
\emph{Comm. Pure Appl. Math.} \textbf{43} (1990), 233--271.
\href{https://doi.org/10.1002/cpa.3160430204}{doi:10.1002/cpa.3160430204}.

\bibitem{LiLu2022}
Y.~Y.~Li and S.~Lu,
Monge--Amp\`ere equation with bounded periodic data,
\emph{Anal. Theory Appl.} \textbf{38} (2022), 128--147.
\href{https://doi.org/10.4208/ata.OA-0022}{doi:10.4208/ata.OA-0022}.

\bibitem{LiNguyenWang2020}
Y.~Y.~Li, L.~Nguyen, and B.~Wang,
Towards a Liouville theorem for continuous viscosity solutions to fully nonlinear elliptic equations in conformal geometry,
in \emph{Geometric Analysis: In Honor of Gang Tian's 60th Birthday},
Progress in Mathematics, vol.~333, Birkh\"auser, Cham, 2020, 221--244.
\href{https://doi.org/10.1007/978-3-030-34953-0_11}{doi:10.1007/978-3-030-34953-0\_11}.

\bibitem{MS1992}
J.~Moser and M.~Struwe,
On a Liouville-type theorem for linear and nonlinear elliptic differential equations on a torus,
\emph{Bol. Soc. Brasil. Mat. (N.S.)} \textbf{23} (1992), 1--20.
\href{https://doi.org/10.1007/BF02584809}{doi:10.1007/BF02584809}.

\bibitem{NguyenII}

N.~C.~Nguyen, On the H\"older continuous subsolution problem for the
complex Monge--Amp\`ere equation, II,
\emph{Anal. PDE} \textbf{13} (2020), 435--453.
\href{https://doi.org/10.2140/apde.2020.13.435}{doi:10.2140/apde.2020.13.435}.


\bibitem{Pogorelov}
A.~V.~Pogorelov,
On the improper convex affine hyperspheres,
\emph{Geom. Dedicata} \textbf{1} (1972), 33--46.
\href{https://doi.org/10.1007/BF00147379}{doi:10.1007/BF00147379}.

\bibitem{RS1984}
D.~Riebesehl and F.~Schulz,
A priori estimates and a Liouville theorem for complex Monge--Amp\`ere equations,
\emph{Math. Z.} \textbf{186} (1984), 57--66.
\href{https://doi.org/10.1007/BF01215491}{doi:10.1007/BF01215491}.

\bibitem{Savin}
O.~Savin, Small perturbation solutions for elliptic equations,
\emph{Comm. Partial Differential Equations} \textbf{32} (2007), 557--578.
\href{https://doi.org/10.1080/03605300500394405}{doi:10.1080/03605300500394405}.

\bibitem{Stein1970}

E.~M.~Stein,
\emph{Singular Integrals and Differentiability Properties of Functions},
Princeton Mathematical Series, vol.~30, Princeton University Press, Princeton, NJ, 1970.
\href{https://doi.org/10.1515/9781400883882}{doi:10.1515/9781400883882}.

\bibitem{Stein1993}
E.~M.~Stein,
\emph{Harmonic Analysis: Real-Variable Methods, Orthogonality, and Oscillatory Integrals},
Princeton Mathematical Series, vol.~43, Princeton University Press, Princeton, NJ, 1993.
\href{https://doi.org/10.1515/9781400883929}{doi:10.1515/9781400883929}.


\bibitem{TZ2016}
E.~V.~Teixeira and L.~Zhang,
Global Monge--Amp\`ere equation with asymptotically periodic data,
\emph{Indiana Univ. Math. J.} \textbf{65} (2016), 399--422.
\href{https://doi.org/10.1512/iumj.2016.65.5687}{doi:10.1512/iumj.2016.65.5687}.

\bibitem{Wang2026}
L.~Wang,
Liouville-type theorems for coupled-drift Monge--Amp\`ere equations,
preprint, 2026.
\href{https://arxiv.org/abs/2608.04478}{arXiv:2608.04478}.

\bibitem{Wang2013}
Y.~Wang,
A Liouville theorem for the complex Monge--Amp\`ere equation,
preprint, 2013.
\href{https://arxiv.org/abs/1303.2403}{arXiv:1303.2403}.

\bibitem{Yau1978}
S.-T.~Yau,
On the Ricci curvature of a compact K\"ahler manifold and the complex Monge--Amp\`ere equation, I,
\emph{Comm. Pure Appl. Math.} \textbf{31} (1978), 339--411.
\href{https://doi.org/10.1002/cpa.3160310304}{doi:10.1002/cpa.3160310304}.

\bibitem{Zeriahi2013}
A.~Zeriahi,
A viscosity approach to degenerate complex Monge--Amp\`ere equations,
\emph{Ann. Fac. Sci. Toulouse Math. (6)} \textbf{22} (2013), 843--913.
\href{https://doi.org/10.5802/afst.1390}{doi:10.5802/afst.1390}.

\bibitem{ZB2018}
W.~Zhang and J.~Bao,
A Calabi theorem for solutions to the parabolic Monge--Amp\`ere equation with periodic data,
\emph{Ann. Inst. H. Poincar\'e C Anal. Non Lin\'eaire} \textbf{35} (2018), 1143--1173.
\href{https://doi.org/10.1016/j.anihpc.2017.09.007}{doi:10.1016/j.anihpc.2017.09.007}.


\bibitem{ZZ2026}
X.~Zhou and R.~Zhu,
Rigidity for convex entire solutions to the $k$-Hessian equation,
preprint, 2026.
\href{https://arxiv.org/abs/2609.20445}{arXiv:2609.20445}.

\end{thebibliography}
\end{document}